\documentclass{gtart_a}
\pdfoutput=1
\usepackage{amssymb,amsmath,amscd,graphicx}
\usepackage{overpic}
\input{diagxy}
\let\to\rightarrow
\let\barrsquare\square
\let\square\undefined

\title{Modifying surfaces in 4--manifolds by twist spinning}

\author{Hee Jung Kim}
\givenname{Hee Jung}
\surname{Kim}
\address{Department of Mathematics\\
McMaster University\\\newline
Hamilton\\
Ontario L8S 4K1\\
Canada}
\email{hjkim@math.mcmaster.ca}
\urladdr{}

\volumenumber{10}
\issuenumber{}
\publicationyear{2006}
\papernumber{2}
\lognumber{0480}
\startpage{27}
\endpage{56}

\doi{}
\MR{}
\Zbl{}

\keyword{Twist spinning}
\keyword{Seiberg--Witten invariants}
\keyword{branched covers}
\keyword{ribbon knots}
\subject{primary}{msc2000}{57R57}
\subject{secondary}{msc2000}{14J80}
\subject{secondary}{msc2000}{57R95}

\received{22 July 2004}
\revised{29 November 2005}
\accepted{2 January 2006}
\published{25 February 2006}
\publishedonline{25 February 2006}
\proposed{Ronald Fintushel}
\seconded{Peter Ozsv\'ath, Ronald Stern}
\corresponding{}
\editor{}
\version{}

\arxivreference{math.GT/0411078}
\arxivpassword{}

\begin{asciiabstract}
In this paper, given a knot K, for any integer m we construct a new
surface Sigma_K(m) from a smoothly embedded surface Sigma in a smooth
4-manifold X by performing a surgery on Sigma. This surgery is based
on a modification of the `rim surgery' which was introduced by
Fintushel and Stern, by doing additional twist spinning. We
investigate the diffeomorphism type and the homeomorphism type of
(X,Sigma) after the surgery. One of the main results is that for
certain pairs (X,Sigma), the smooth type of Sigma_K(m) can be easily
distinguished by the Alexander polynomial of the knot K and the
homeomorphism type depends on the number of twist and the knot. In
particular, we get new examples of knotted surfaces in CP^2, not
isotopic to complex curves, but which are topologically unknotted.
\\
In this paper, given a knot $K$, for any integer $m$ we construct a
new surface $\Sigma_K(m)$ from a smoothly embedded surface $\Sigma$
in a smooth 4--manifold $X$ by performing a surgery on $\Sigma$. This
surgery is based on a modification of the `rim surgery' which was
introduced by Fintushel and Stern, by doing additional twist
spinning. We investigate the diffeomorphism type and the
homeomorphism type of $(X,\Sigma)$ after the surgery. One of the
main results is that for certain pairs $(X,\Sigma)$, the smooth type
of $\Sigma_K(m)$ can be easily distinguished by the Alexander
polynomial of the knot $K$ and the homeomorphism type depends on the
number of twist and the knot. In particular, we get new examples of
knotted surfaces in $\mathbb{C}\mathbf{P}^2$, not isotopic to complex
curves, but which are topologically unknotted.
\end{asciiabstract}

\begin{htmlabstract}
In this paper, given a knot K, for any integer m we construct a
new surface &Sigma;<sub>K</sub>(m) from a smoothly embedded surface &Sigma;
in a smooth 4&ndash;manifold X by performing a surgery on &Sigma;. This
surgery is based on a modification of the `rim surgery' which was
introduced by Fintushel and Stern, by doing additional twist
spinning. We investigate the diffeomorphism type and the
homeomorphism type of (X,&Sigma;) after the surgery. One of the
main results is that for certain pairs (X,&Sigma;), the smooth type
of &Sigma;<sub>K</sub>(m) can be easily distinguished by the Alexander
polynomial of the knot K and the homeomorphism type depends on the
number of twist and the knot. In particular, we get new examples of
knotted surfaces in <b>C</b><b>P</b><sup>2</sup>, not isotopic to complex
curves, but which are topologically unknotted.
\end{htmlabstract}

\AtBeginDocument{\def\notin{\not\in}}

\newtheorem{Thm}{Theorem}[section]
\newtheorem{Prop}{Proposition}[section]
\newtheorem{Lemma}{Lemma}[section]
\newtheorem{Cor}{Corollary}[section]

\theoremstyle{remark}
\newtheorem{Def}{Definition}[section]
\newtheorem{Ex}{Example}[section]
\newtheorem*{Rmk}{Remark}

\makeatletter
\let\c@Prop\c@Thm
\let\c@Lemma\c@Thm
\let\c@Cor\c@Thm
\let\c@Guess\c@Thm
\let\c@Def\c@Thm
\let\c@Ex\c@Thm
\makeatother
\makeautorefname{Thm}{Theorem}
\makeautorefname{Prop}{Proposition}
\makeautorefname{Lemma}{Lemma}
\makeautorefname{Cor}{Corollary}
\makeautorefname{Guess}{Guess}
\makeautorefname{Def}{Definition}
\makeautorefname{Ex}{Example}

\makeop{im}
\makeop{Im}
\makeop{cl}
\makeop{Hom}
\makeop{id}
\makeop{spin}
\makeop{Wh}

\begin{document}

\begin{abstract}
In this paper, given a knot $K$, for any integer $m$ we construct a
new surface $\Sigma_K(m)$ from a smoothly embedded surface $\Sigma$
in a smooth 4--manifold $X$ by performing a surgery on $\Sigma$. This
surgery is based on a modification of the `rim surgery' which was
introduced by Fintushel and Stern, by doing additional twist
spinning. We investigate the diffeomorphism type and the
homeomorphism type of $(X,\Sigma)$ after the surgery. One of the
main results is that for certain pairs $(X,\Sigma)$, the smooth type
of $\Sigma_K(m)$ can be easily distinguished by the Alexander
polynomial of the knot $K$ and the homeomorphism type depends on the
number of twist and the knot. In particular, we get new examples of
knotted surfaces in $\mathbb{C}\mathbf{P}^2$, not isotopic to complex
curves, but which are topologically unknotted.
\end{abstract}

\maketitle

\section{Introduction}
\label{sec:1}

Let $X$ be a smooth 4--manifold and $\Sigma$ be an embedded positive
genus surface and nonnegative self-intersection. In \cite{FS1},
Fintushel and Stern introduced a technique, called `rim surgery', of
modifying $\Sigma$ without changing the ambient space $X$. This
surgery on $\Sigma$ may change the diffeomorphism type of the
embedding $\Sigma_K$ but the topological embedding is preserved when
$\pi_1(X-\Sigma)$ is trivial. Rim surgery is determined by a knotted
arc $K_+\in B^3$, and may be described as follows. Choose a curve
$\alpha$ in $\Sigma$, which has a neighborhood $S^1\times B^3$
meeting $\Sigma$ on an annulus $S^1\times I$. Replacing the pair
$(S^1\times B^3,S^1\times I)$ by $(S^1\times B^3,S^1\times K_+)$
gives a new surface $\Sigma_K$ in $X$.

In \cite{Z}, Zeeman described the process of twist-spinning an
$n$--knot to obtain an $(n{+}1)$--knot. Here an {\em $n$--knot} is a
locally flat pair $(S^{n+2}, K)$ with $K\cong S^n$. Then here is the
description for the process of twist-spinning to obtain a knot in
dimension 4: Suppose we have a knotted arc $K_+$ in the half 3--space
$\mathbb{R}_{+}^3$, with its end points in $\mathbb{R}^2=\partial\mathbb{R}_{+}^3$. Spinning $\mathbb{R}_{+}^3$ about
$\mathbb{R}^2$ generates $\mathbb{R}^4$, the arc $K_+$ generates a
knotted 2--sphere in $\mathbb{R}^4$, called a {\em spun knot}. During
the spinning process we spin the arc $K_+$ $m$ times keeping its end
points within $\mathbb{R}_{+}^3$, obtaining again a 2--sphere $K(m)$
in $\mathbb{R}^4$. A more explicit definition is the following.

For any 1--knot $(S^3, K)$, let $(B^3, K_{+})$ be its ball pair
with the knotted arc $K_{+}$. Let $\tau$ be the diffeomorphism of
$(B^3, K_{+})$, called `twist map' defined in \fullref{sec:2}. Then for
any integer $m$ this induces a 2--knot called the $m$--\emph{twist spun
knot}
$$(S^4,K(m))=\partial(B^3, K_{+})\times B^2\cup_{\partial}
  (B^3,K_{+})\times_{\tau^m}\partial B^2$$
where $(B^3,K_{+})\times_{\tau^m}\partial B^2$ means that
$(B^3,K_{+})\times [0,1]/(x,0)=(\tau^m x,1)$.

In this paper, using these two ideas --- rim surgery and spun knot --- we
will construct a new surface, denoted by $\Sigma_K(m)$, from the embedded
surface in $X$ without changing its ambient space.  Our technique may be
called a `twist rim surgery'. We will see later (in \fullref{sec:3} and
\fullref{sec:4}) that the smooth and topological type of $\Sigma_K(m)$
obtained by twist rim surgery depends on $m$, $K$, and $\Sigma$. For
a precise definition of the surgery, we will give two descriptions of
$\Sigma_K(m)$. One is provided by using the twist map $\tau$ in the
construction of Zeeman's twist spun knot. The other one can be obtained
by performing the same operation which Fintushel and Stern introduced
in \cite{FS3} as it corresponds to doing a surgery on a homologically
essential torus in $X$. In \cite{FS3}, they constructed exotic manifolds
$X_{K}$ according to a knot $K$ and also showed that the Alexander
polynomial $\Delta_K(t)$ of $K$ can detect the smooth type of $X_K$.

In our circumstance, we consider a pair $(X,\Sigma)$, where $X$ is a
smooth simply connected $4$--manifold and $\Sigma$ is an embedded
genus $g$ surface with self-intersection $n\geq{0}$ such that the
homology class $[\Sigma]=d\cdot{\beta}$, where $\beta$ is a
primitive element in $H_{2}(X)$ and $\pi_1(X-\Sigma)=\mathbb{Z}/d$.
Then in \fullref{sec:3}, we will study the smooth type of $\Sigma_K(m)$
obtained by performing twist rim surgery on $\Sigma$. In fact, using
the result in \cite{FS1}, we conclude that the Alexander polynomial
$\Delta_K(t)$ of $K$ can distinguish the smooth type of
$\Sigma_K(m)$.
  In particular, applying this result to $\mathbb{C}\mathbf{P}^2$ we can get new
examples of knotted surfaces in $\mathbb{C}\mathbf{P}^2$, not isotopic to complex
curves. This solves, for an algebraic curve of degree $\ge 3$,
Problem 4.110 in the Kirby list \cite{K}. Note that $d=1,2$
which are the only degrees where the curve is a sphere, are still open.

In \fullref{sec:4}, we will study topological conditions under which
 $(X,\Sigma_K(m))$ is pairwise homeomorphic to $(X,\Sigma)$.
This problem is also related to the knot type of $K$ and the
relation between $d$ and $m$. In particular, if $d\not\equiv\pm 1
\pmod{m}$ then computing the fundamental group of the exterior of
surfaces in $X$ we easily distinguish $(X,\Sigma_K(m))$ and
$(X,\Sigma)$ for some nontrivial knot $K$. But when $d\equiv\pm 1
\pmod{m}$, it turns out that the fundamental group
$\pi_1(X-\Sigma_K(m))$ is same as $\pi_1(X-\Sigma)=\mathbb{Z}/d$.
So, in the case $d\equiv\pm 1 \pmod{m}$ we show that if $K$ is a
ribbon knot and the $d$--fold cover of the knot complement $S^3-K$ is
a homology circle then $(X,\Sigma)$ and $(X,\Sigma_K(m))$ are
topologically equivalent. This means that there is a pairwise
homeomorphism
$(X,\Sigma)\longrightarrow (X,\Sigma_K(m))$.

\section{Definitions}
\label{sec:2}

Let $X$ be a smooth 4--manifold and let $\Sigma$ be an embedded
surface of positive genus $g$. Given a knot $K$ in $S^3$, let $E(K)$
be the exterior $\cl(S^3-K\times D^2)$ of $K$. First we need to
consider a certain diffeomorphism $\tau$ on $(S^3,K)$ which will be
used to define our surgery. Take a tubular neighborhood of the knot
and then using a suitable trivialization with 0--framing, let
$\partial E(K)\times I=K\times{\partial{D^2}}\times{I}$ be a collar
of $\partial E(K)$ in $E(K)$ with $\partial E(K)$ identified with
$\partial E(K)\times\{0\}$. Define $\tau\co (S^3,K)\longrightarrow
(S^3,K)$ by
\begin{equation}\label{tau}
\tau(x\times{e^{i\theta}}\times{t}) =x\times{e^{i(\theta + 2\pi
t)}}\times{t} \quad \mbox{for} \quad
x\times{e^{i\theta}}\times{t}\in K\times{\partial{D^2}}\times{I}
\end{equation}
and $\tau(y)=y$ for $y\notin K\times{\partial{D^2}}\times{I}$.

Note that $\tau$ is not the identity on the collar $\partial
E(K)\times I=K\times{\partial{D^2}}\times{I}$. However, it is the
identity on the exterior $\cl(S^3-K\times{\partial{D^2}}\times{I})$
of the collar. If we restrict $\tau$ to the exterior of the knot $K$
then $\tau$ is isotopic to the identity although the isotopy is not
the identity on the boundary of the knot complement. Explicitly, the
isotopy can be given as the following. For any $s\in [0,1]$,
$$\tau_s(x\times{e^{i\theta}}\times{t})=x\times{e^{i\theta+2\pi t(1-s)+2\pi s}}
\times{t}.$$
We will refer to this diffeomorphism $\tau$ as
a {\em twist map}.

Now take a non-separating curve $\alpha$ in $\Sigma$. Then choose
a trivialization of the normal bundle $\nu(\Sigma)|_{\alpha}$ in
$X$, $\alpha\times I\times D^2=\alpha\times B^3\longrightarrow
\nu(\Sigma)|_{\alpha}$ where $\alpha\times I$ corresponds to the
normal bundle $\nu(\alpha)$ in $\Sigma$. For any trivialization of
the tubular neighborhood of $\alpha$ we can construct a new
surface from $\Sigma$ using the chosen curve. We will choose a
specific framing of $\alpha$ later in \fullref{sec:3} to study the
diffeomorphism type of the new surface constructed in the way
discussed now. Identifying $\alpha$ with $S^1$, two descriptions
of the construction of $(X,\Sigma_K(m))$ called
$m$--\emph{twist rim surgery} follow.

\begin{Def}\label{def:first}
Define for any integer $m$,
$$(X,\Sigma_K(m))=(X,\Sigma) -
S^1\times(B^3,I)\cup_{\partial} S^1\times_{\tau^m} (B^3,K_{+}).$$
\end{Def}

Note that for $m=0$, $\Sigma_K(m)$ is the surface obtained by rim
surgery. In \cite{FS1}, its smooth type was studied when
$\pi_1(X-\Sigma)$ is trivial. As in the paper \cite{FS1}, we will
consider the smooth type of the new surface obtained by $m$--twist
surgery in the extended case
where $\pi_1(X-\Sigma)$ is cyclic.

If $\alpha$ is a trivial curve, that is it bounds a disk in
$\Sigma$, we can simply write $(X,\Sigma_K(m))$ as the following.

\begin{Lemma}
If $\alpha$ is a trivial curve
in $\Sigma$, then $(X,\Sigma_K(m))$ is the connected sum
$(X,\Sigma)$ with the $m$--twist spun knot $(S^4, K(m))$ of
$(S^3,K)$.
\end{Lemma}

\begin{proof}
Considering the decomposition of $(X,\Sigma_K(m))$ in \fullref{def:first}.
$$(X,\Sigma_K(m))=(X,\Sigma) -
  S^1\times(B^3,I)\cup_{\partial} S^1\times_{\tau^m} (B^3,K_{+}),$$
we write the boundary of the ball $(B^3, I)$ in the definition as
$$\partial (B^3,I)=(S^2, \{N,S\})=(D_{+}^2,\{N\})\cup (D_{-}^2,\{S\})$$
where $D_{+}^2$, $D_{-}^2$ are 2--disks and N, S are north and south
poles respectively. Also recall that we identified $\alpha$ as $S^1$
in the definition and by the choice of $\alpha$, let's denote the
disk bounded by $\alpha$ as $B^2$ in $\Sigma$. Then we can rewrite
\begin{multline*}
(X,\Sigma_K(m)) =\\
  \bigl((X,\Sigma){-}(S^1{\times}(B^3,I) \cup
  B^2{\times} (D_{+}^2,\{N\}))\bigr) \cup
  \bigl(B^2{\times}(D_{+}^2,\{N\})\cup S^1{\times_{\tau^m}} (B^3,K_{+})\bigr).
\end{multline*}
Note that the first component of this decomposition is
$$(X,\Sigma)-S^1\times (B^3,I)\cup_{\partial B^2\times D_{+}^2}
  B^2\times (D_{+}^2,\{N\})=(X,\Sigma)-(B^4, B^2).$$
In the second component
$$B^2\times (D_{+}^2,\{N\})\cup_{\partial B^2\times D_{+}^2}
S^1\times_{\tau^m} (B^3,K_{+}),$$
gluing $B^2\times
(D_{-}^2,\{S\})$ to $B^2\times (D_{+}^2,\{N\})$ along $B^2\times
\partial D_{+}^2$ and then taking it out later again we can write
\begin{multline*}
\bigl(B^2 {\times} (D_{+}^2,\{N\})\bigr)
  \cup_{B^2{\times} \partial D_{+}^2}
\bigl(B^2{\times}
(D_{-}^2,\{S\})\bigr) \cup_{\partial}
  \bigl(S^1{\times}_{\tau^m} (B^3\!\!,\!K_{+})\bigr){-}
  \bigl(B^2{\times} (D_{-}^2,\{S\})\bigr)\\
=\bigl(B^2{\times}\partial(B^3\!\!,\!K_{+})\bigl)\cup_{\partial}
  \bigl(S^1{\times}_{\tau^m}(B^3\!\!,\!K_{+})\bigr){-}
  \bigl(B^2{\times} (D_{-}^2,\{S\})\bigr).
\end{multline*}
Considering the definition of twist spun knot in \fullref{sec:1} we can
realize this is
$$\bigl(S^4, K(m)\bigr)-\bigl(B^2{\times} (D_{-}^2,\{S\})\bigr).$$
So,
$$(X,\Sigma_K(m))=\bigl((X,\Sigma) -(B^4,B^2)\bigr)\cup
  \bigl((S^4, K(m))-B^2{\times}(D_{-}^2,\{S\})\bigr)$$
where the union is taken along the boundary.
\end{proof}

Let's move on to another description of $(X,\Sigma_K(m))$ which is
useful in distinguishing the diffeomorphism types of $\Sigma_K(m)$.
For a non-separating curve $\alpha$ in $\Sigma$, after a
trivialization, the normal bundle $\alpha$ in $X$ is of the form
$\alpha\times I\times D^2=\alpha\times B^3$ where $\alpha\times I$
in $\Sigma$. Consider $\alpha\times \gamma\subset \alpha\times
I\times D^2$ where $\gamma$ is a pushed-in copy of the meridian
circle $\{0\}\times\partial D^2\subset I\times D^2$. Under our
trivialization, $\alpha\times \gamma$ is diffeomorphic to a torus
$T$ in $X-\Sigma$, called a {\em rim torus} by Fintushel and Stern.
Note that this torus $T$ is nullhomologous in $X$. Let $N(\gamma)$
be a tubular neighborhood of $\gamma$ in $B^3=I\times D^2$ and
$\gamma'$ be the curve $\gamma$ pushed off into $\partial
N(\gamma)$. Then we will identify $\alpha\times N(\gamma)$ as a
neighborhood $N(T)$ of $T$ under the trivialization so that
$\alpha\times
N(\gamma)\subset\nu(\Sigma)|_{\alpha}\subset\nu(\Sigma)$. For a knot
$K$ in $S^3$, let's denote by $\mu_{K}$ the meridian and
$\lambda_{K}$ the longitude of the knot. Now consider the following
manifold
$$\alpha\times (B^3-N(\gamma))\cup_{\varphi}(S^1\times E(K))$$
where the gluing map $\varphi$ is the diffeomorphism determined by
$\varphi_{*}(\alpha)= m\mu_{K}+S^1$,
$\varphi_{*}(\gamma')=\mu_{K}$, and $\varphi_{*}(\partial D^2)=\lambda_{K}$.

\begin{Def}\label{def:second}
Suppose that $T\cong \alpha\times\gamma$ is the smooth torus in
$X$ as above. Define
$$(X,\Sigma_K(m))=(X - N(T),\Sigma)\cup_{\varphi}(E(K)\times S^1,\emptyset).$$
\end{Def}

This description means that performing a surgery on a smooth torus
$T$ in $X$, we obtain $X$ again but $\Sigma$ might be changed. Now
we need to check those two descriptions are the same definitions
for our construction.

\begin{Lemma}\label{lem:equivdef}
\fullref{def:first} and \fullref{def:second} are equivalent.
\end{Lemma}

\begin{proof} Given a knot $K$, recall that knotting the arc
$I=I\times \{0\}\subset B^3=I\times B^2$ can be achieved by a
cut-paste operation on the complement. Let $\gamma$ be an unknot
which is the meridian of the arc $I$ in $B^3$, $E(K)$ be the
exterior of the knot $K$ in $S^3$ and $N(\gamma)$ be the tubular
neighborhood of $\gamma$ in $B^3$. If we replace the tubular
neighborhood $N(\gamma)$ by $E(K)$ then we get $B^3$ with the
knotted arc $K_+$ instead of the trivial arc $I$. More precisely,
note that $(B^3,K_{+})=(\nu(\partial{B^3}\cup{K_{+}}),K_{+})\cup
E(K)$ where $\nu(\partial{B^3}\cup{K_{+}})$ is the normal bundle in
$B^3$ (see \fullref{fig:1}). Let $\gamma'$ be the push off of
$\gamma$ onto $\partial N(\gamma)$.

\begin{figure}[ht!]
\begin{center}
\begin{overpic}[scale=0.7]{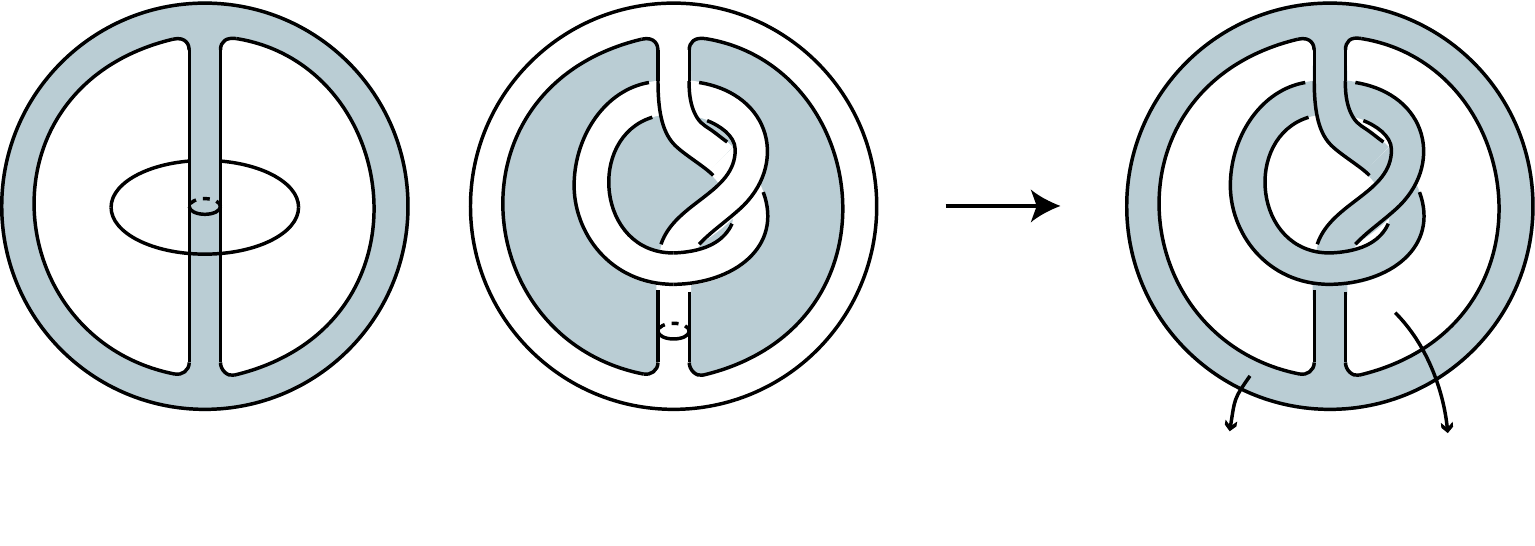}
\put(1,4){\small $(B^3,I)-N(\gamma)$}
\put(14.5,20.5){\small$\gamma'$}
\put(20,20.5){\small $\gamma$}
\put(28,4){\small $\cup_f$}
\put(38,4){\small $E(K)$}
\put(45,13){\small $\mu_k$}
\put(71,4){\small $v(\partial B^3 \cup K_-) \cup E(K)$}
\end{overpic}
\caption{Diffeomorphism $h\co (B^3-N(\gamma), I)\cup_{f} E(K)\to (B^3,K_{+})$}
\label{fig:1}
\end{center}
\end{figure}

Then there is a diffeomorphism $(B^3-N(\gamma), I)\rightarrow
(\nu(\partial{B^3}\cup{K_{+}}),K_{+})$ mapping $\gamma'$ to
$\mu_{K}$ which induces a diffeomorphism
$$h\co  (B^3-N(\gamma), I)\cup_{f} E(K)\longrightarrow
(\nu(\partial{B^3}\cup{K_{+}}),K_{+})\cup E(K)=(B^3,K_{+}),$$
where $f\co \partial N(\gamma)\longrightarrow \partial E(K)$ is a
diffeomorphism determined by identifying $\gamma'$ to ${\mu_{K}}$.
Note that the diffeomorphism $h$ has $h(I)=K_{+}$ and $h|_{E(K)}=\id$.

Recalling the map $\tau$ defined in \eqref{tau}, we note that $h$ is
the identity on $E(K)$ but $\tau$ is not, whereas on the outside of
$E(K)$, $\tau$ is the identity but $h$ is not. This implies that
$\tau$ is equivariant with respect to $h$, $\tau\circ
h$=$h\circ\tau$. This $\tau$ induces a well-defined diffeomorphism
mapping $[x,t]$ to $[h(x),t]$
$$\bigl(((B^3,I)-N(\gamma)) \cup_{f}
  E(K)\bigr)\times_{\tau^m}S^1\longrightarrow
  (B^3,K_{+})\times_{\tau^m}S^1.$$
Since $\tau^m$ is the identity on $(B^3,I)-N(\gamma)$,
$(((B^3,I)-N(\gamma)) \cup_{f} E(K))\times_{\tau^m}S^1$ is the same
as $((B^3,I)-N(\gamma))\times S^1 \cup_{f\times
1_{S^1}}(E(K)\times_{\tau^m}S^1)$ and thus we have
$$((B^3,I)-N(\gamma))\times S^1 \cup_{f\times
  1_{S^1}}(E(K)\times_{\tau^m}S^1)\longrightarrow
  (B^3,K_{+})\times_{\tau^m}S^1.$$
Extending by the identity gives a diffeomorphism
\begin{multline*}
((X,\Sigma) - (B^3,I)\times S^1)\cup_{\partial}
((B^3,I)-N(\gamma))\times S^1
\cup_{f\times 1_{S^1}}(E(K)\times_{\tau^m}S^1)\longrightarrow \\
((X,\Sigma)-(B^3,I)\times S^1)
\cup_{\partial}(B^3,K_+)\times_{\tau^m} S^1.
\end{multline*}
Rewriting
\begin{equation*}
\begin{split}
((X,\Sigma) -(B^3,I)&\times S^1)\cup_{\partial}
((B^3,I)-N(\gamma))\times S^1
\cup_{f\times 1_{S^1}}(E(K)\times_{\tau^m}S^1)\\
&=X-N(\gamma)\times S^1
\cup_{f\times 1_{S^1}}(E(K)\times_{\tau^m}S^1)\\
&=X-\gamma\times D^2\times S^1 \cup_{f\times
1_{S^1}}(E(K)\times_{\tau^m}S^1),
\end{split}
\end{equation*}
we get a diffeomorphism
$$X{-}\gamma{\times} D^2\times S^1
\cup_{f{\times} 1_{S^1}}(E(K){\times_{\tau^m}}S^1)\rightarrow
((X,\Sigma){-}(B^3\!,\!I)\times S^1)
\cup_{\partial}(B^3\!,\!K_+){\times_{\tau^m}}S^1.$$
Note that here the gluing map $f\times 1_{S^1}$ sends $\alpha$ to
$S^1$, $\gamma'$ to  $\mu_{K}$ and $\partial{D^2}$ to
${\lambda_{K}}$ where $\mu_{K}$ and $\lambda_{K}$ are the meridian
and the longitude of the knot $K$. Since $\tau^m$ is isotopic to
identity, the isotopy induces a diffeomorphism $E(K)\times
S^1\longrightarrow E(K)\times_{\tau^m}S^1$. Again extending by the
identity gives a diffeomorphism
$$X-\gamma \times D^2\times S^1 \cup_{f\times
1_{S^1}}(E(K)\times_{\tau^m}S^1)
\rightarrow
(X - \gamma\times D^2\times S^1)\cup_{\varphi}(E(K)\times S^1),$$
where $\varphi$ is given by
\begin{align*}
\alpha &\longleftrightarrow{S^1+m\mu_{K}}\\
\gamma'&\longleftrightarrow{\mu_{K}}\\
\partial{D^2}&\longleftrightarrow{\lambda_{K}}.
\end{align*}
Therefore the result follows.
\end{proof}

\section{Diffeomorphism types}
\label{sec:3}

Now let $X$ be a smooth simply connected $4$--manifold and $\Sigma$
an embedded genus $g$ surface with self-intersection $n\geq{0}$ and
homology class $[\Sigma]=d\cdot{\beta}$, where $\beta$ is a
primitive element in $H_{2}(X)$ and $\pi_1(X-\Sigma)=\mathbb{Z}/d$.
Since $\Sigma$ is diffeomorphic to $T^2\#\cdot\cdot\cdot\# T^2$,
let's choose a curve $\alpha$ whose image is the curve $\{pt\}
\times{S^1}$ in the first $T^2=S^1\times S^1$. As we discussed in
the previous section, a neighborhood of $\alpha$ in $X$ is of the
form $\alpha\times I\times D^2=\alpha\times B^3$, where
$\alpha\times I$ is in $\Sigma$. But we need to choose a certain
trivialization of the normal bundle $\nu(\alpha\times I)$ in $X$
which will be used in \fullref{sec:4} when we compute some topological
invariants to identify the homeomorphism type of $\Sigma_K(m)$. It is possible to choose a trivialization $\sigma$ of
$\nu(\alpha\times I)$ with the property that for some point $p\in
\partial D^2$, $\sigma|\alpha\times\{0\}\times p$ is trivial in
$H_1(X-\Sigma)$; we arbitrarily choose one trivialization
$\sigma\co \alpha\times I\times D^2\longrightarrow \nu(\alpha\times I)$
and let $\alpha'$ be $\sigma|\alpha\times \{0\}\times p$ for some
$p\in \partial D^2$. By composing $\sigma$ with a self
diffeomorphism of $\alpha\times I\times D^2$ sending the element
$(e^{i\theta},t,z)$ to $(e^{i\theta},t,e^{ik\theta}z)$ for an
appropriate integer $k$, we can arrange $\alpha'$ to be the zero
homology element in $H_1(X-\Sigma)\cong\mathbb{Z}/d$, that is
generated by the meridian
$\sigma(pt\times\partial D^2 )$ of $\Sigma$.

For a given $d$, the relation between $\Sigma_K(m)$ and $\Sigma$
depends somewhat on $m$. For example, if $d\not\equiv\pm 1 \pmod{m}$
then for a nontrivial knot $K$, the surface $\Sigma_K(m)$ can
be distinguished (even up to homeomorphism) from $\Sigma$ by
considering the fundamental group $\pi_1(X-\Sigma_K(m))$.  First, we need
to understand the explicit expression of this group.

In this paper, we will denote by $(X,Y)^d$ a $d$--fold covering of
$X$ branched along $Y$.
\begin{Lemma} \label{lem:FG}
Let $\mu$ be the meridian of the knotted arc $K_+$ and let the
base point $*$ be in $\partial
E(K)=K\times\partial{D^2}\times\{0\}$. Then
$$\pi_1( X{-}\Sigma_{K}(m)) = \bigl\langle\pi_1(B^3{-}K_+,*) ~|~
\mu^{d}=1,\beta =\tau_{*}^{m}(\beta), \text{for all } \beta\in
\pi_1(B^3{-}K_+,*)\bigr\rangle.$$
\end{Lemma}

\begin{proof} Considering the definition of
$(X,\Sigma_{K}(m))$,
  we have that the complement of $\Sigma_{K}(m)$ in $X$,
$X-\Sigma_{K}(m)$, is $(X-S^1\times
B^3-\Sigma)\cup{S^1\times_{\tau^m} (B^3-K_{+})}$. Then we get that
the intersection of the two components in the decomposition is
$$(X-S^1\times B^3-\Sigma)\cap{S^1\times_{\tau^m} (B^3-K_{+})}=
S^1\times(\partial{B^3}-\{\text{two points}\}).$$ Here we need to note
that the action of $\tau$ on $\partial{B^3}-\{$ two points $\}$ is
trivial. Then using Van Kampen's theorem for this decomposition, we
have the following diagram:
$$\bfig
  \barrsquare<1500,500>[
  \pi_1(S^1{\times}(\partial B^3 - \{\text{two points}\}))`
  \pi_1(X-S^1\times B^3-\Sigma)`
  \pi_1(S^1\times_{\tau^m} (B^3-K_{+}))`
  \pi_1(X - \Sigma_{K}(m));
  \varphi_1`\varphi_2`\psi_1`\psi_2]
  \efig$$
Note that $X{-}S^1{\times}B^3{-}\Sigma$ is homotopy equivalent to
$X{-}\Sigma$ and $\pi_1(X{-}\Sigma)\cong\mathbb{Z}/d$ is
generated by the meridian $\gamma$ of $\Sigma$.  We also know that
$\pi_1 (S^1{\times}(\partial{B^3}{-}\{\text{two points}\}))$ is generated
by $[S^1]$, which is identified with the class of the curve
$\alpha'$ pushed off along a given trivialization of neighborhood of
$\alpha$, and by $\mu$. Since the meridian $\mu$ of the knot is
identified with $\gamma$, $\varphi_1$ is onto and so $\psi_{2}$ is
also onto. Moreover,
  $\ker \psi_{2}=\langle\varphi_2(\ker\varphi_{1})\rangle$.
Since $\ker\varphi_{1}=\langle\alpha' , \mu^d\rangle$ and
\begin{multline*}
\pi_1(S^1{\times_{\tau^m}}
(B^3,K_{+}))=\\
\bigl\langle\pi_1(B^3{-}K_+), \alpha'  ~|~
{\alpha'}^{-1}\beta\alpha'=\tau_{*}^{m}(\beta) \text{ for all }\beta \in
\pi_1(B^3{-}K_+)\bigr\rangle,
\end{multline*}
it follows that
\begin{align*}
&\pi_1(X-\Sigma_{K}(m))\\
&=\bigl\langle\pi_1(B^3{-}K_+), \alpha'  ~|~\alpha'=1 , \mu^d=1,
  {\alpha'}^{-1}\beta\alpha'=\tau_{*}^{m}(\beta)  \text{ for all }
  \beta \in \pi_1(B^3{-}K_+)\bigr\rangle\\
&= \bigl\langle\pi_1(B^3{-}K_+) ~|~  \mu^d =1, \beta =\tau_{*}^{m}(\beta)
\text{ for all } \beta \in \pi_1(B^3{-}K_+)\bigr\rangle.
\end{align*}
which completes the proof.
\end{proof}

The following example shows that we can distinguish $\Sigma_K(m)$
using
$\pi_1$.

\begin{Ex} For any nontrivial knot $K$, let $d=2$,
ie $\pi_1(X-\Sigma)=\mathbb{Z}/2$, and let $m$ be any even number.
If we consider the fundamental group $\pi_1(X-\Sigma_K(m))$, then by
\fullref{lem:FG},
$$\pi_1(X{-}\Sigma_{K}(m))=\bigl\langle\pi_1(B^3{-}K_+,*)
  ~|~ \mu^{d}=1,\beta =\tau_{*}^{m}(\beta)
  \text{ for all }\beta\in \pi_1(B^3{-}K_+,*)\bigr\rangle,$$
where $\mu$ is the meridian of the knotted arc $K_+$ and the base
point $*$ is in
$\partial E(K)=K\times\partial{D^2}\times\{0\}$.

Recall the group of the knot $\pi_1(B^3-K_+,*)$ has the Wirtinger
presentation
$$\langle g_1,g_2,\ldots ,g_n \mid r_1,r_2,\ldots,r_n\rangle,$$
where $g_1=\mu$ and other generators $g_i$ represent the loop
that, starting from a base point, goes straight to the $i^{th}$
over-passing arc in
the knot diagram, encircles it and returns to the base point.

Note that $\tau_{*}^{m}(g_1)=g_1$ and
$\tau_{*}^{m}(g_i)=g_1^{-m}g_ig_1^{m}$ for other generators $g_i$ by
the definition of $\tau$. Since $d=2$ ie $g_1^2=1$ and $m$ is an
even number, $\tau_{*}^{m}(g_i)=g_1^{-m}g_ig_1^{m}$ is always $g_i$
and thus we get
$$\pi_1(X-\Sigma_{K}(m))=\pi_1(B^3-K_+)/\mu^2=\pi_1(S^3-K)/\mu^2.$$
If we take a 2--fold branched cover $(S^3,K)^2$ along the knot $K$
then the fundamental group $\pi_1((S^3,K)^2)$ is same as the group
$\pi_1((S^3-K)^2)/\wtilde \mu$, where $(S^3-K)^2$ is the 2--fold
unbranched cover and $\wtilde \mu$ is a lift of $\mu$. So
$\pi_1(S^3-K)/\mu^2$ has $\pi_1((S^3,K)^2)$ as an index 2 subgroup.
The Smith conjecture \cite{MB} states that for any $d\ge 1$, the
fundamental group of a $d$--fold branched cover $\pi_1((S^3,K)^d)$ is
nontrivial unless $K$ is a trivial knot. Hence
$\pi_1(X-\Sigma_{K}(m))$ has a nontrivial index 2 subgroup and so
$\pi_1(X-\Sigma_{K}(m))\not\cong\mathbb{Z}/2$. This proves that
there is no homeomorphism
$(X-\Sigma)\to (X-\Sigma_{K}(m))$.
\end{Ex}

A more interesting  case is when $\pi_1$ does not distinguish the
embedding of $\Sigma_K(m)$, so that we have to use other means to
show that $\Sigma$ is not diffeomorphic to $\Sigma_K(m)$. In
particular, for the case $d\equiv\pm 1\pmod{m}$, we have:

\begin{Prop} \label{prop:FG}
If $d\equiv\pm 1 \pmod{m}$ then
$\pi_1(X-\Sigma)=\pi_1(X-\Sigma_{K}(m))=\mathbb{Z}/d$.
\end{Prop}

\begin{proof}
If $d=1$ then by \fullref{lem:FG},
$\pi_1(X-\Sigma)=\pi_1(X-\Sigma_{K}(m))=\{1\}$. So, we assume $d>1$.
To express $\pi_1(X-\Sigma)$ more explicitly, in a Wirtinger
presentation of the knot group $\pi_1(B^3-K_+,*)$, choose meridians
$g_j$ conjugate to the meridian $g_1=\mu$ of the knot $K$ for each
$j=2,...,n$ as generators of $\pi_1(B^3-K_+) $. Then with
\fullref{lem:FG}, we represent $\pi_1(X-\Sigma_K(m))$ by
$$\langle g_1,g_2,\ldots ,g_n \mid g_{1}^d=1 , r_1,\ldots ,r_n, \beta
=\tau_{*}^{m}(\beta) \text{ for all }\beta\in \pi_1(B^3-K_+)\rangle$$
where $r_1,\ldots,r_n$ are relations of $\pi_1(B^3-K_+)$.

Considering the definition of $\tau$, $\tau_{*}(\mu)=\mu$ and
$\tau_{*}(g_j)=\mu^{-1}g_j\mu$ for each $j=2,\ldots,n$ so that we
rewrite
\begin{equation*}
\begin{split}
\pi_1(&X-\Sigma_{K}(m)) \\ & =\ \langle g_1,g_2,\ldots ,g_n  \mid
  g_{1}^d=1, r_1,\ldots ,r_n, g_j=g_{1}^{-m}g_j g_{1}^m
\mbox{\ for\ } j=2,\ldots ,n\rangle.
\end{split}
\end{equation*}
Now we claim that this is equal to $\langle g_1,g_2,\ldots ,g_n  \mid
g_{1}^d, r_1,\ldots ,r_n,
g_1=g_{1}^{-1}g_j g_{1}$  for $j=2,\ldots ,n\rangle$.

Since $d\equiv\pm 1 \pmod{m}$, we can write $d=mk\pm 1$ for some
integer $k$. Let $l=d-m$. Then $l=d-m=mk\pm 1-m=m(k-1)\pm 1$.
\begin{align*}
g_j  =  g_{1}^{-m}g_j g_{1}^m &\Longrightarrow
g_{1}^{-l}g_jg_{1}^{l}=g_{1}^{-l}(g_{1}^{-m}g_j g_{1}^m) g_{1}^{l}\\
&\Longrightarrow g_{1}^{-l}g_j g_{1}^{l}=g_{1}^{-(l+m)}g_j
g_{1}^{(l+m)}=g_j
&(l+m=d)\\
&\Longrightarrow g_{1}^{-m(k-1)\mp 1}g_j g_{1}^{m(k-1)\pm 1}=g_j
&(l=m(k-1)\pm 1)\\
&\Longrightarrow g_{1}^{\mp 1}(g_{1}^{-m(k-1)}g_j g_{1}^{m(k-1)})
g_{1}^{\pm 1}=g_j\ldots &(*)
\end{align*}
We claim that  $g_{1}^{-m(k-1)}g_j g_{1}^{m(k-1)}=g_j$; if $k-1=0$
or $1$ then it is clearly true. Let's assume that it is true for
$k-1=i$. For $k-1=i+1$, by induction
$$g_{1}^{-m(i+1)}g_j g_{1}^{m(i+1)}=g_{1}^{-mi}(g_{1}^{-m}g_jg_{1}^{m})
g_{1}^{mi}=g_{1}^{-mi}g_j g_{1}^{mi}=g_j.$$
This implies that $(*)$ becomes $g_{1}^{\mp 1}g_j g_{1}^{\pm 1}=g_j$
and so we now get
\begin{equation*}
\begin{split}
\pi_1(&X-\Sigma_{K}(m))\\
&=\ \langle g_1,g_2,\dots ,g_n  \mid  g_{1}^d =1 , r_1,\ldots ,r_n,
[g_1,g_j]=1  \mbox{\ for\ } j=2,\ldots ,n\rangle.
\end{split}
\end{equation*}
If we consider the Wirtinger presentation of the knot group then we
can show $g_1=g_2=...=g_n$ with the relations $r_1,..,r_n$ and
$[g_1,g_j]$; corresponding to the following crossing, the relator
gives $g_2 g_s=g_s g_1$ or $g_s g_2=g_1 g_s$.

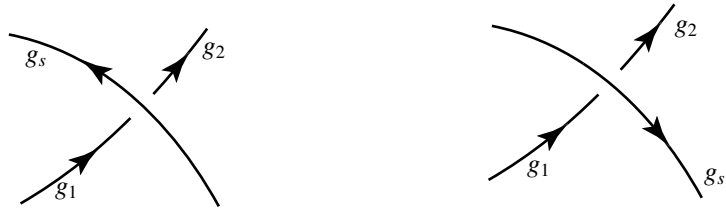
\begin{figure}[ht!]
\begin{center}
\begin{overpic}{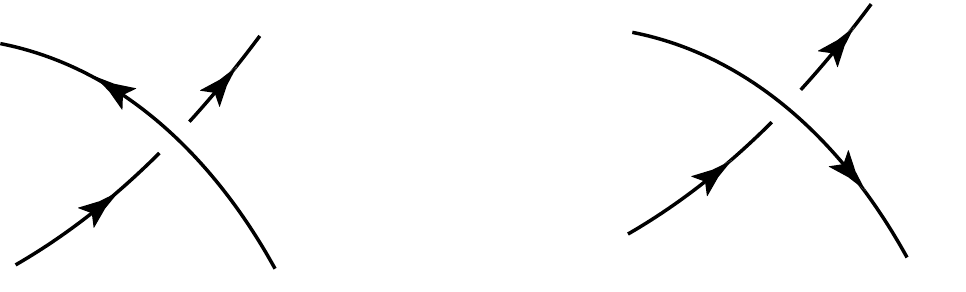}
\put(2,22){\small $g_s$}
\put(6,4){\small $g_1$}
\put(26,23){\small $g_2$}
\put(70,7){\small $g_1$}
\put(94,6){\small $g_s$}
\put(90,26){\small $g_2$}
\end{overpic}
\caption{Wirtinger presentation of the knot group}
\end{center}
\end{figure}

So, $g_1=g_2$. By an induction argument, we can conclude that
$g_1=g_2=\ldots =g_n$. This proves that
$$\pi_1(X-\Sigma_{K}(m))\ =\ \langle\mu  \mid  \mu^d =1 \rangle  \ \cong \
\mathbb{Z}/d.\proved$$
\end{proof}

\begin{Rmk}  The same technique works for many other
cases, for example if $d=2$
and $m$ is an odd integer. 
\end{Rmk}

We can also distinguish some $\Sigma_{K}(m)$ smoothly by using
relative Seiberg--Witten (SW) theory, following the technique of
Fintushel and Stern \cite{FS2}. In \cite{FS3}, they introduced a
method called `knot surgery' modifying a 4--manifold while preserving
its homotopy type by using a knot in $S^3$ and also gave a formula
for the SW-invariant of the new manifold to detect the
diffeomorphism type under suitable circumstances.

Let $X$ be a smooth 4--manifold and $T$ in $X$ be an imbedded 2--torus
with trivial normal bundle. (In \cite{T3}, C Taubes showed the
`$c$--embedded' condition on the torus in the original paper \cite{FS3}
to be unnecessary.) Then the knot surgery may be described as
follows.

Let $K$ be a knot in $S^3$, and $K\times D^2$ be the trivialization
of its open tubular neighborhood given by the 0--framing. Let
$\varphi\co \partial(T\times D^2)\longrightarrow \partial(K\times
D^2)\times S^1$ be any  diffeomorphism with $\varphi(p\times
\partial D^2)=K\times q$ where $p\in T$, $q\in\partial D^2\times
S^1$ are points. Define
$$X_K=(X-T\times D^2)\cup_{\varphi}E(K)\times S^1.$$
In our situation, the surgical construction of $\Sigma_K(m)$ is
performing a surgery on a torus $T$ in $X$ called a `rim torus'.
Recall the torus $T$ has the form $\gamma\times\alpha$ where
$\gamma$ is the meridian of $\Sigma$ and $\alpha$ is a curve in
$\Sigma$ (see \fullref{lem:equivdef}). In other words, we remove a
neighborhood of the torus and sew in $E(K)\times S^1$ along the
gluing map given in \fullref{def:second}. Considering this
identification, we can
observe that the pair $(X,\Sigma_K(m))$ is obtained by a knot surgery.

Fintushel and Stern wrote a note to fill a gap in the proof of the
main theorem in \cite{FS1}. In the note \cite{FS2}, they explained
the effect of rim surgery on the relative Seiberg--Witten invariant
of $X-\Sigma$. The $m$--twist rim surgery on $X-\Sigma$ affects its
relative Seiberg--Witten invariant exactly same as rim surgery. So we
will refer to the note \cite{FS2} to distinguish the pairs
$(X,\Sigma)$ and $(X,\Sigma_K(m))$ smoothly.

If the self-intersection $\Sigma\cdot\Sigma =n\geq{0}$, blow up $X$
$n$ times to get a pair $(X_n,\Sigma_n)$ and reduce the self
intersection to zero. For simplicity, we may assume that
$\Sigma\cdot\Sigma =0$. In general, the relative Seiberg--Witten
invariant $SW_{X,\Sigma}$ is an element in the Floer homology of the
boundary $\Sigma\times S^1$ \cite{KM}. We restrict $SW_{X,\Sigma}$
to the set $\mathcal{T}$ which is the collection of
$\spin^{c}$--structures $\tau$ on $X-N(\Sigma)$ whose restriction to
$\partial{N(\Sigma)}$ is the $\spin^c$--structure $\pm{s_{g-1}}$
corresponding to the element $(g-1,0)$ of $H^2(\Sigma\times
S^1)\cong \mathbb{Z}\oplus H^1(\Sigma)$. Then we obtain a
well-defined integer-valued Seiberg--Witten invariant
$SW^{\mathcal{T}}_{X,\Sigma}$ and so get a Laurent polynomial
$SW^{\mathcal{T}}_{X,\Sigma}$ with variables in
$$A=\{\alpha\in H^2(X-\Sigma) | \alpha|_{\Sigma\times
S^1}=\pm{s_{g-1}}\}.$$
If there is a diffeomorphism $f\co (X,\Sigma)\to
(X',\Sigma')$ then it induces a map $f^*\co A'\to A$ sending
$SW^{\mathcal{T}}_{X',\Sigma'}$ to $SW^{\mathcal{T}}_{X,\Sigma}$.

\begin{Thm} \label{thm:main}
Suppose the relative Seiberg--Witten invariant
$SW^{\mathcal{T}}_{X,\Sigma}$ is nontrivial. If there is a
diffeomorphism $(X,\Sigma_{K}(m))\longrightarrow (X,\Sigma_{J}(m))$
then the set of coefficients (with multiplicity) of $\Delta_{K}(t)$
is equal to that of $\Delta_{J}(t)$, where $\Delta_{K}(t)$ and
$\Delta_{J}(t)$ are the Alexander polynomials of $K$ and $J$
respectively.
\end{Thm}

\begin{proof} If there is a pairwise diffeomorphism
$(X,\Sigma_{K}(m))\longrightarrow (X,\Sigma_{J}(m))$ then it induces
a diffeomorphism $(X_n,\Sigma_{n,K}(m))\longrightarrow
(X_n,\Sigma_{n,J}(m))$. So, we now may assume that
$\Sigma\cdot\Sigma=0$.

According to the note \cite{FS2}, the proof of the knot surgery
theorem \cite{FS3} works in the relative case to show that
$$SW^{\mathcal{T}}_{(X-\Sigma)_{K}}
  =SW^{\mathcal{T}}_{X,\Sigma}\cdot \Delta_{K}(r^2)$$
where $r=[T]$ is the element of $R$, the subgroup of $H^2(X-\Sigma)$
generated by the rim torus $T$ of $\Sigma$. Note that the rim torus
$T$ is homologically essential in $X-\Sigma$.

Since the relative Seiberg--Witten invariant
$SW^{\mathcal{T}}_{X,\Sigma_{K}(m)}=SW^{\mathcal{T}}_{(X-\Sigma)_{K}}$,
applying the knot surgery theorem to the $m$--twist rim surgery we also
get that the coefficients of $SW^{\mathcal{T}}_{X,\Sigma}\cdot
\Delta_{K}(r^2)$ must be equal to those of
$SW^{\mathcal{T}}_{X,\Sigma}\cdot \Delta_{J}(r'^2)$.
\end{proof}

\begin{Rmk}
\begin{enumerate}
\item The theorem implies that for
$\Delta_{K}(t)\neq 1$, $(X,\Sigma)$ is not pairwise
diffeomorphic to $(X,\Sigma_{K}(m))$.

\item In \cite{FS1} standard pairs $(Y_g,S_g)$ were defined where $Y_g$
is a simply connected K\"ahler surface, $S_g$ is a primitively
embedded genus $g\geq 1$ Riemann surface in $Y_g$ with $S_g\cdot
S_g=0$. According to the note \cite{FS2}, the hypothesis
$SW_{X\#_{\Sigma=S_g}Y_g}\ne 1$ of \cite{FS1} implies
$SW^{\mathcal{T}}_{X,\Sigma}\ne 1$ by the gluing formula \cite{KM}.

\item $SW_{X\#_{\Sigma=S_g}Y_g}$ is nontrivial when $\Sigma$ is a
complex curve in a complex surface.
\end{enumerate}
\end{Rmk}

The case of curves in $\mathbb{C}\mathbf{P}^2$ is particularly interesting. By
applying \fullref{thm:main}, we obtain the following corollary.

\begin{Cor}
For $d>2$ with $d\equiv\pm 1 \pmod{m}$, if $\Sigma$ is a degree
$d$--curve in $\mathbb{C}\mathbf{P}^2$ then $(\mathbb{C}\mathbf{P}^2,\Sigma)$ is not
pairwise diffeomorphic to $(\mathbb{C}\mathbf{P}^2,\Sigma_{K}(m))$ for any
knot $K$ with $\Delta_K(t)\neq{1}$, but
$\pi_1(\mathbb{C}\mathbf{P}^2-\Sigma_K(m))\cong\mathbb{Z}/d$.
\end{Cor}

\begin{proof}
Note that $\Sigma$ is a symplectically embedded surface with positive
genus $g=\frac{1}{2}(d-1)(d-2)$. Under the construction
in~\cite{FS1}, $S_g$ is also symplectically embedded in $Y_g$ since
$S_g$ is a complex submanifold of the K\"ahler manifold $Y_g$. Since
the group $\pi_1(C\mathbf{P}^2-\Sigma)=\mathbb{Z}/d$, note that
$\pi_1(C\mathbf{P}^2-\Sigma_{K}(m))=\mathbb{Z}/d$ by
\fullref{prop:FG}.

Let us denote by $CP^{2}_{d^2}$ the manifold obtained by blowing up
$d^2$ times $CP^{2}$. Then $CP^{2}_{d^2}\#_{\Sigma_{d^2}=S_g}Y_g$
is also a symplectic manifold by Gompf \cite{G}. So (see Taubes~\cite{T1}),
$$SW_{CP^{2}_{d^2}\#_{\Sigma_{d^2}=S_g}Y_g}\ne{0}.$$
By \fullref{thm:main}, the result follows.
\end{proof}

This means that for any $d\geq 3$, there are infinitely many smooth
oriented closed surfaces $\Sigma$ in $\mathbb{C}\mathbf{P}^2$ representing the
class $dh\in H_2(\mathbb{C}\mathbf{P}^2)$, where $h$ is a generator of
$H_2(\mathbb{C}\mathbf{P}^2)$, having
$\mbox{genus}(\Sigma)=\frac{1}{2}(d-1)(d-2)$ and
$\pi_1(\mathbb{C}\mathbf{P}^2-\Sigma)\cong\mathbb{Z}/d$, such that the pairs
$(\mathbb{C}\mathbf{P}^2,\Sigma)$ are pairwise smoothly non-equivalent. Such
examples, for $d\ge 5$, were known by the work of Finashin which we
describe in order to contrast it with our construction. In
\cite{Finashin}, he constructed a new surface by knotting a standard
one along a suitable annulus
membrane.

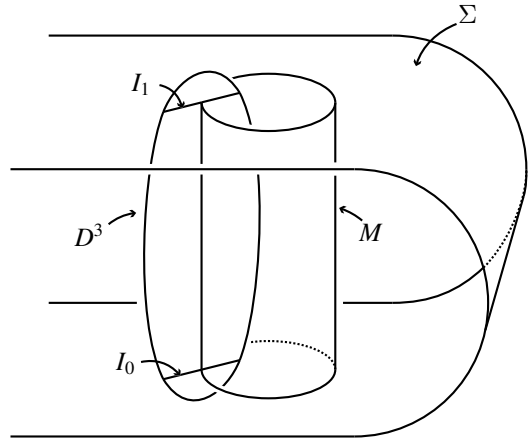
\begin{figure}[ht!]
\begin{center}
\begin{overpic}[scale=0.8]{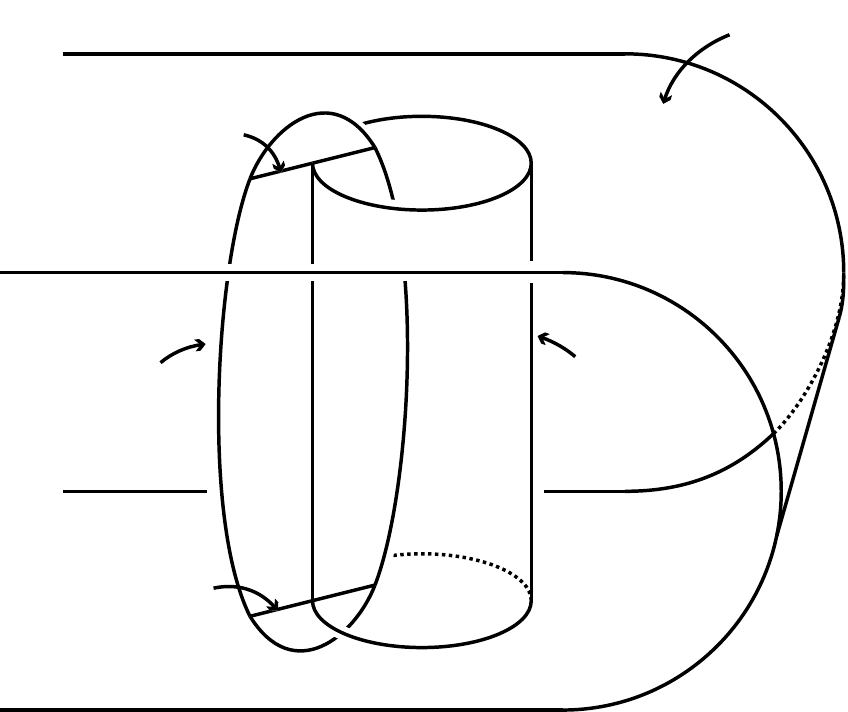}
\put(12,37){\small $D^3$}
\put(23,67){\small $I_1$}
\put(20,13){\small $I_0$}
\put(67,38){\small $M$}
\put(86,80){\small $\Sigma$}
\end{overpic}
\caption{$(\mathbb{C}\mathbf{P}^2,\Sigma)$}\label{fig:fina}
\end{center}
\end{figure}

More precisely, let $X$ be a 4--manifold and $\Sigma$ be a smoothly
embedded surface. Suppose that there is a smoothly embedded surface
$M$ in $X$, called a `membrane', such that $M\cong S^1\times I$,
$M\cap \Sigma=\partial M$ and $M$ meets to $\Sigma$ normally along
$\partial M$. By adjusting a trivialization of its regular
neighborhood $U$, we can assume that $U(\cong S^1\times D^3)\cap
\Sigma=S^1\times f$, where
  $f=I_0\sqcup I_1=I\times\partial I$ is a disjoint union of
two unknotted segments of a part of the boundary  of a band
$b=I\times I$ in $D^3$. Here the band $b=I\times I$ is trivially
embedded in $D^3$ and the intersection $I\times I\cap \partial
D^3=\partial I\times I$ (see \fullref{fig:fina}).

\begin{figure}[ht!]
\begin{center}
\begin{overpic}[scale=.8]{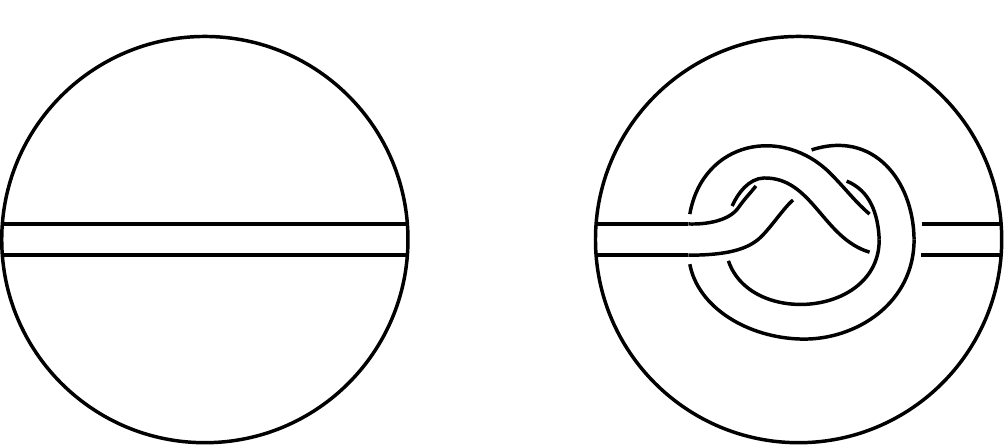}
\put(8,40){\small $D^3$}
\put(68,40){\small $D^3$}
\put(18,23.5){\small $I_1$}
\put(18,15){\small $I_0$}
\end{overpic}
\caption{$(D^3,I\times I)$ and $(D^3,K_{+}\times I)$} \label{fig:ballpair}
\end{center}
\end{figure}
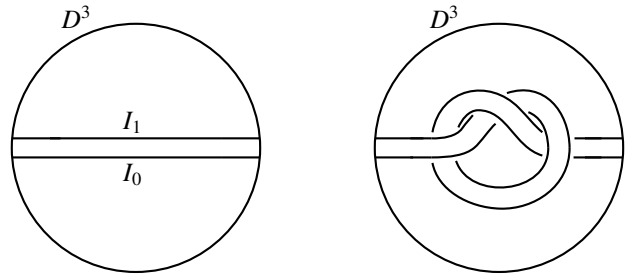

Then given a knot $K$ in $S^3$, we can get a new surface
$\Sigma_{K,F}$ by knotting $f$ along $K$ in $D^3$ (see \fullref{fig:ballpair}).

In \cite{Finashin}, Finashin showed that we can find such a membrane
$M$ in $\mathbb{C}\mathbf{P}^2$ and proved that  $(\mathbb{C}\mathbf{P}^2,
\Sigma_{K,F})$ is pairwise non-equivalent to
$(\mathbb{C}\mathbf{P}^2,\Sigma)$ for an algebraic curve $\Sigma$ of
degree $d\ge 5$.
In particular, for an even degree he showed that the double cover
branched along $\Sigma_{K,F}$ is diffeomorphic to the 4--manifold
obtained from the double cover branched along $\Sigma$ by knot
surgery along the torus, which is the pre-image of the membrane $M$
in the covering, via the knot $K\#K$. So, the knot surgery theorem
in \cite{FS3} distinguishes the branch covers by comparing their
SW-invariants. For odd cases, one can use the same argument using
$d$--fold coverings
to show smooth non-equivalence of embeddings.

Our examples constructed by twist spinning are different from
Finashin's for a degree $d\ge 5$. To see this, we compute the
SW-invariant of the branched cover of
$(\mathbb{C}\mathbf{P}^2,\Sigma_{K}(m))$. Let $Y$ be a $d$--fold branch cover
along $\Sigma$ and $Y_{K,m}$ be a $d$--fold branch cover along
$\Sigma_{K}(m)$. Let's consider the description for the branch
cover $Y_{K,m}$. We write $Y_{K,m}$ as the union of two $d$--fold
branched covers:
$$(Y_{K,m},\Sigma_K(m))=(X-S^1\times B^3,\Sigma-S^1\times I)^d\cup_{\partial}
  (S^1\times_{\tau^m}(B^3,K_{+}))^d$$
Since the homology group $H_1(X-S^1\times B^3-\Sigma)\cong
H_1(X-\Sigma)\cong\mathbb{Z}/d$, the branch cover $(X-S^1\times
B^3,\Sigma-S^1\times I)^d$ is unique and is the same as $Y-S^1\times
B^3$. We also need to note that
$(S^1\times_{\tau^m}(B^3,K_{+}))^d=S^1\times_{\wtilde\tau^m}(B^3,K_{+})^d
$ for some lift $\wtilde\tau^m$ of $\tau^m$ which is referred to in
the proof for \fullref{prop:cobo}. So we rewrite
$$(Y_{K,m},\Sigma_K(m))=((Y,\Sigma)-S^1\times (B^3,I))\cup_{S^1\times
  S^2}(S^1\times_{\wtilde\tau^m}(B^3,K_{+})^d).$$
If $K$ is any knot with the homology $H_1((S^3-K)^d)\cong\mathbb{Z}$
then $S^1\times_{\wtilde\tau^m}(B^3,K_{+})^d$ is homologically
equivalent to $S^1\times B^3$. We may look at knots, introduced in
\fullref{sec:4}, having the property that their $d$--fold covers are homology
circles. An extension of the result of Vidussi in \cite{V} shows
$$SW_{Y_{K,m}}=SW_{Y}.$$
But the SW-invariant of branched cover along the surface
$\Sigma_{K,F}$ constructed by Finashin is not standard as we saw
above.  Our examples also cover the case of degree $d=3$ and
$4$ which were not treated in his paper.

\begin{Rmk}
By the same argument in Fintushel and Stern
\cite{FS1}, we can also say that if $X$ is a simply connected
symplectic 4--manifold and $\Sigma$ is a symplectically embedded
surface then $\Sigma_K(m)$ is not smoothly ambient isotopic to a
symplectic submanifold of $X$ for $\Delta_K(t)\neq 1$. Using Taubes'
result in \cite{T1}, we can easily get a proof of this (see
\cite{FS1} for more detail).
\end{Rmk}

\section{Homeomorphism types}
\label{sec:4}

In this section, we shall investigate when $\Sigma_{K}(m)$ is
topologically equivalent to $\Sigma$. As we saw in the previous
section, in the case $d\equiv\pm 1 \pmod{m}$ their complements in
$X$ have the same fundamental group. So, for this case one would
like to show that they are pairwise homeomorphic under a certain
condition by constructing an explicit $s$--cobordism. Note that it is
not known if Finashin's examples are topologically unknotted
\cite[Remark, p50]{Finashin}. Recall that the $s$--cobordism theorem
gives a way for
showing manifolds are homeomorphic.

Let $W$ be a compact $n$--manifold with the boundary being the
disjoint union of manifolds $M_0$ and $M_1$. Then the original
$s$--cobordism theorem states that for $n\geq 6$, $W$ is diffeomorphic
to $M_0\times[0,1]$ exactly when the inclusions of $M_0$ and $M_1$
in $W$ are homotopy equivalences and the Whitehead torsion
$\tau(W,M_0)$ in $\Wh(\pi_1(W))$ is zero. By the work of M Freedman
\cite{Freedman}, the $s$--cobordism theorem is known to hold
topologically in the case $n=5$ when $\pi_1(W)$ is
poly-(finite or cyclic). A relative $s$--cobordism theorem also holds.

To make use of those theorems we shall construct a relative
$h$--cobordism from $X-\nu(\Sigma)$ to $X-\nu(\Sigma_{K}(m))$ and then
apply
the relative $s$--cobordism theorem.

\begin{figure}[ht!]
\begin{center}
\begin{overpic}[scale=0.7]{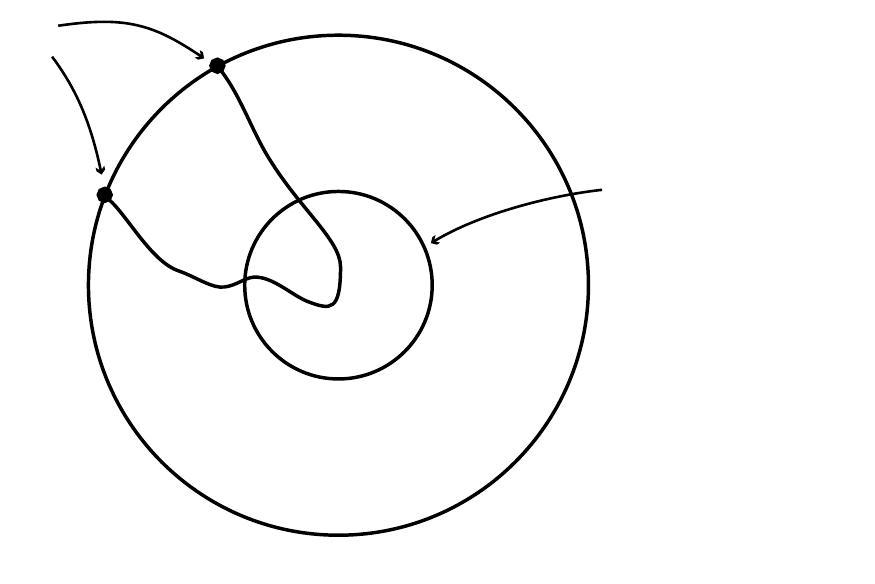}
\put(1,62){\small $K$}
\put(70,43){\small $(B',B'\cap\Delta)$}
\put(60,10){\small $(B^4,\Delta)$}
\end{overpic}
\caption{Ribbon disk in $B^4$} \label{fig:ball}
\end{center}
\end{figure}

First consider the following situation. Let $K$ be a ribbon knot in
$S^3$ so that $(S^3,K)=\partial(B^4,\Delta)$ for some ribbon disc
$\Delta$ in $B^4$. By Lemma 3.1 in \cite{GR},
$\pi_1(S^3-K)\longrightarrow\pi_1(B^4-\Delta)$ is surjective. Take
out a 4--ball $(B', B'\cap\Delta)$ from the interior of
$(B^4,\Delta)$ such that $B'\cap\Delta$ is an unknotted disk (see
\fullref{fig:ball}).

Let $A=\Delta-(B'\cap\Delta)$ then we can easily note that $A$ is a
concordance between $K$ and an unknot $O$. Let $K=K_+\cup K_-$ where
$K_+$ is a knotted arc and $K_{-}$ is a trivial arc diffeomorphic to
$I$. Write $S^3=B_{+}^3\cup B_{-}^3$ where $B_{+}^3$, $B_{-}^3$ are
3--balls. Let's assume that $B_{-}^3\times I\subset S^3\times I$ with
$(B_{-}^3\times I,B_{-}^3\times I\cap A)=(B_{-}^3\times I,I\times
I)$ and
$(B_{-}^3\times1,B_{-}^3\times1\cap A)=(B_{-}^3\times1,K_{-})$.

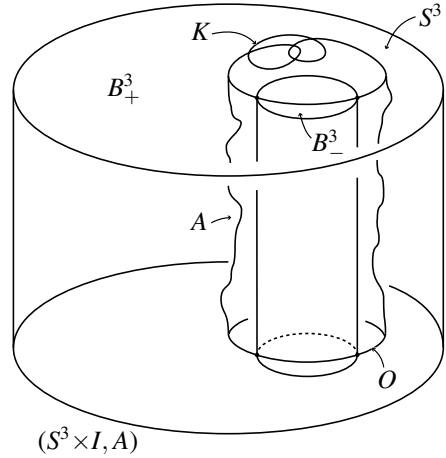
\begin{figure}[ht!]
\begin{center}
\begin{overpic}[scale=.6]{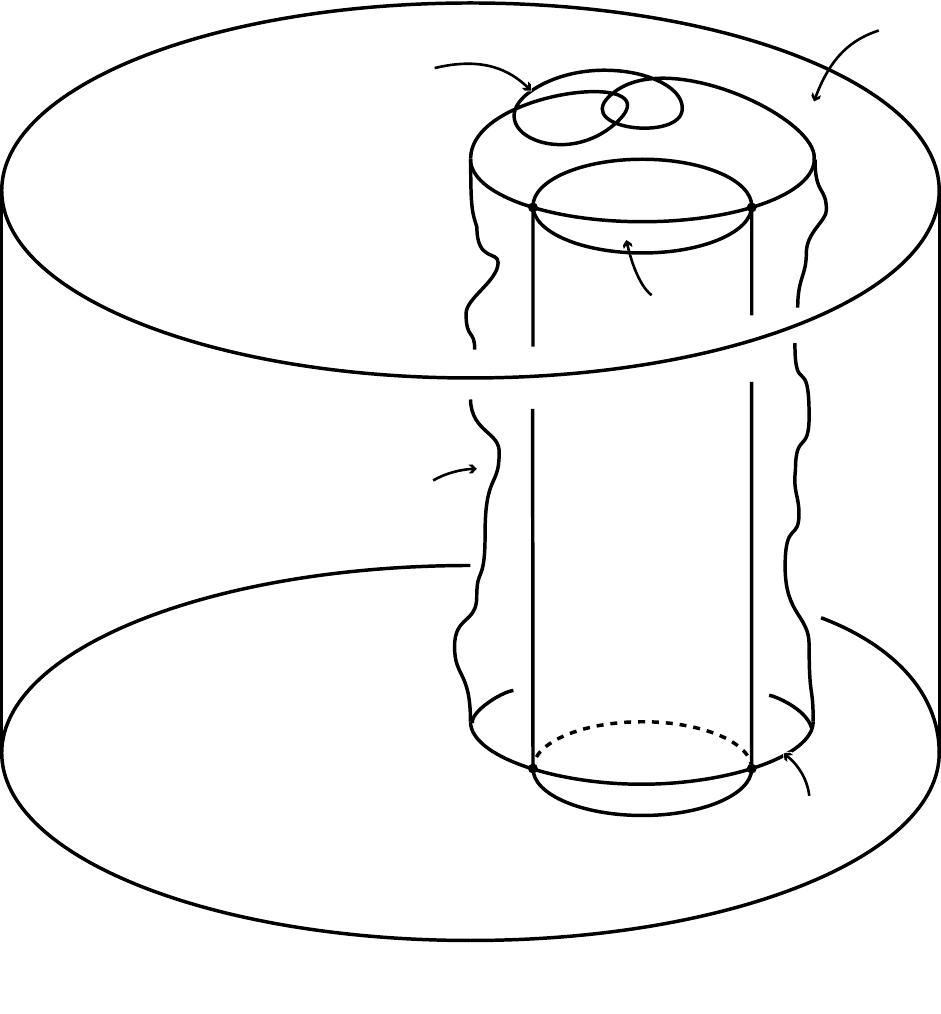}
\put(20,80){\small $B^3_+$}
\put(64,68){\small $B^3_-$}
\put(5,4){\small $(S^3{\times}I,A)$}
\put(78,17){\small $O$}
\put(38,92){\small $K$}
\put(87,95){\small $S^3$}
\put(38,51){\small $A$}
\end{overpic}
\caption{A concordance between $K$ and unknot}\label{fig:cobo}
\end{center}
\end{figure}

If we take out $B_{-}^3\times I$ from $S^3\times I$ then we are left with
  $(S^3\times I,A)-(B_{-}^3\times I,I\times I)=
(B_{+}^3\times I, A-I\times I)$.
  Denoting  $A-I\times I$ by $A_+$, we have
$B_{+}^3\times 1\cap A_{+}=K_{+}$ and $B_{+}^3\times 0\cap
A_{+}=O_{+}$ where $O_+$ is a trivial arc of $O$ (see
\fullref{fig:cobo}).

We will define a self diffeomorphism on $(S^3\times I,A)$ in the
same way that we defined the {\em twist map} in \fullref{sec:2}. Recall
$(S^3\times I,A)-(B_{-}^3\times I,I\times I)= (B_{+}^3\times I,
A_{+})$. Note the normal bundle $\nu(A)$ in $S^3\times I$ is
$A\times D^2$ and let $E(A)$ be the exterior $\cl(S^3\times I-A\times
D^2)$ of $A$ in $S^3\times I$. Then $E(A)$ coincides (up to
isotopy), with $\cl(B_{+}^3\times I-A_{+}\times D^2)$. Thus,
$\partial E(A)=A\times
\partial D^2$ is
$\partial(\cl(B_{+}^3\times I-A_{+}\times D^2))\cong T\times I$
where $T$ is a torus. Let $A\times \partial D^2\times I$ be the
collar of
  $\partial E(A)$ in $E(A)$. Define
$\tau\co (S^3\times I,A)\longrightarrow (S^3\times I,A)$ by
$$\tau(x\times{e^{i\theta}}\times{t}) =x\times{e^{i(\theta + 2\pi
  t)}}\times{t} \quad \mbox{for} \quad
  x\times{e^{i\theta}}\times{t}\in A\times{\partial{D^2}}\times{I}$$
and $\tau(y)=y$ for $y\notin A\times{\partial{D^2}}\times{I}$.

Then note that $\tau$ is the identity on a neighborhood of $A_+$ and that
$\tau|_{B_{+}^3\times 0-O_{+}}$ and $\tau|_{B_{+}^3\times 1-K_{+}}$ are
the twist maps induced by the unknot $O$ and the knot $K$ defined in
\eqref{tau}. Denote those maps by  $\tau_O$ and $\tau_K$
respectively.
Using this diffeomorphism $\tau$, we can also
construct a new submanifold $(\Sigma\times I)_{A}(m)$ from an
embedded manifold $\Sigma\times I$ to
$X\times I$ in the way to construct a new surface $\Sigma_K(m)$.

\begin{Def}\label{def:cobo} Under the above notation, define
$$(X\times I, (\Sigma\times I)_{A}(m))=X\times I-S^1\times(B^3\times
  I, I\times I) \cup S^1\times_{\tau^m}(B^3\times I, A_{+}).$$
\end{Def}

\begin{figure}[ht!]
\begin{center}
\begin{overpic}[scale=0.7]{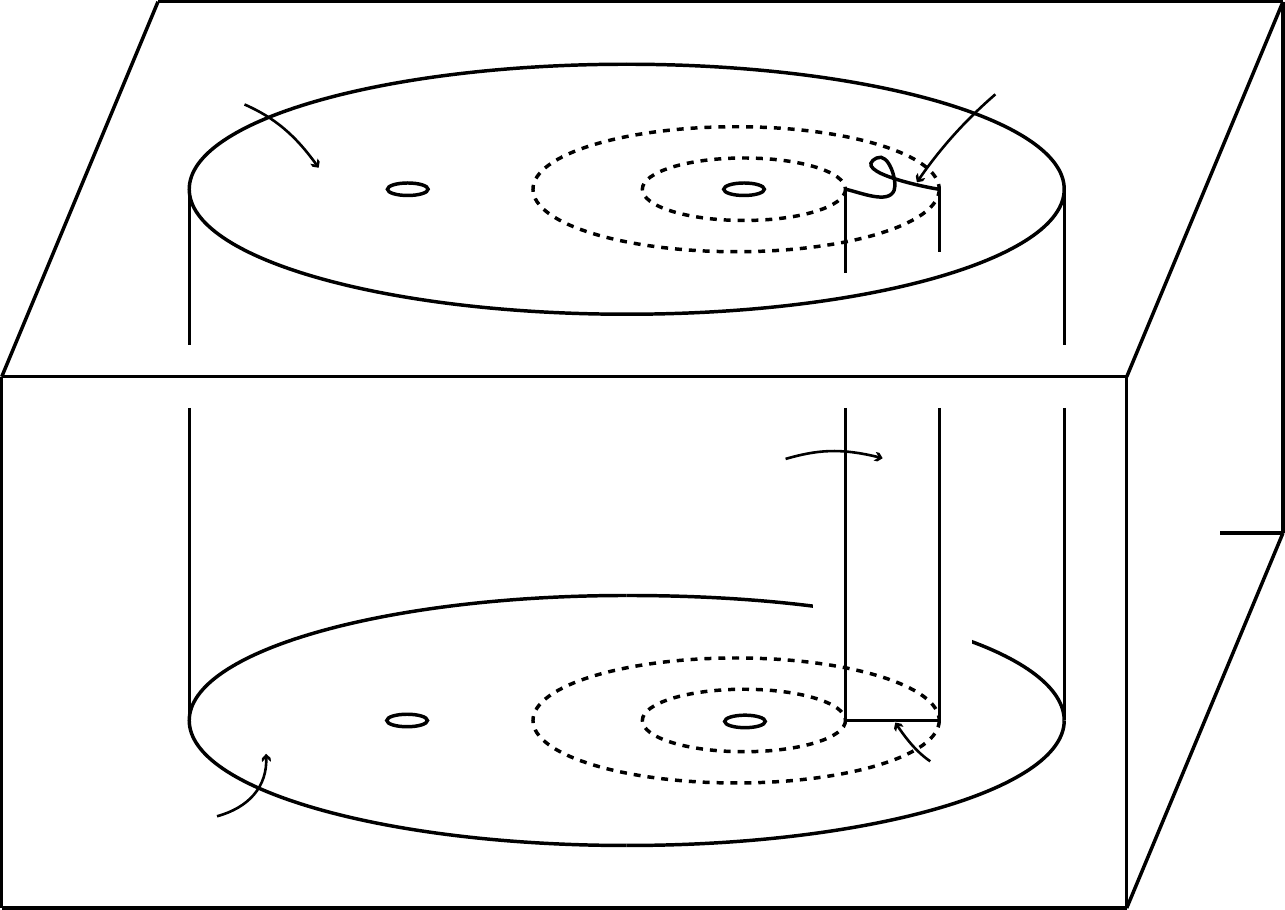}
\put(14,64){\small $\Sigma_k(m)$}
\put(13,5){\small $\Sigma$}
\put(57,34){\small $A_-$}
\put(84,2){\small $X$}
\put(84,43){\small $X$}
\put(77,64){\small $K_+$}
\put(72,10){\small $O_-$}
\end{overpic}
\caption{A cobordism between  $(X,\Sigma)$ and $(X,\Sigma_{K}(m))$}\label{Fig:cobor}
\end{center}
\end{figure}
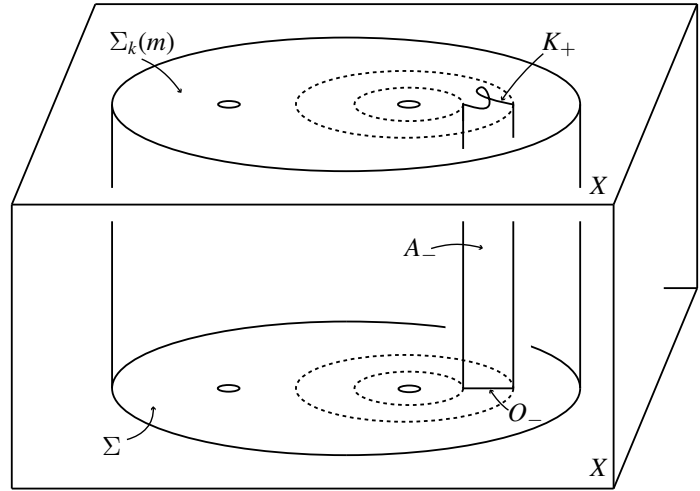

Then we can easily note that
\begin{align*}
X\times 1 &= X-S^1\times(B^3\times 1, I\times 1)
  \cup S^1\times_{\tau_{K}^m}(B^3\times 1,
  K_{+})=(X,\Sigma_{K}(m)), \\
X\times 0 &= X-S^1\times(B^3\times 0, I\times 0)
  \cup S^1\times_{\tau_{O}^m}(B^3\times 0, O_{+})=(X,\Sigma)
\end{align*}
and so the complement $X\times I-(\Sigma\times I)_{A}(m)$ gives a
concordance between $X-\Sigma$ and $X-\Sigma_{K}(m)$ (See
\fullref{Fig:cobor}). We will denote this concordance by $W$ and
will later show this $W$ is a $h$--cobordism under certain conditions.
Here we note that the cobordism $W$ is a product near the boundary.
To see what conditions are needed, consider several other properties
first.

Recall for any pair $(X,Y)$, we denote by $X^d$ a $d$--fold cover of
$X$ and $(X,Y)^d$ a $d$--fold cover of $X$ branched along $Y$. We
know $H_*(S^3-K)\to H_*(B^4-\Delta)$ is an isomorphism
but generally, $H_*((S^3-K)^d)\to H_*((B^4-\Delta)^d)$
is not.  It is true
when $K$ is a ribbon knot:

\begin{Lemma} \label{lem:HC}
If $K$ is a ribbon knot and the homology of $d$--fold cover of
$S^3-K$, $H_1((S^3-K)^d)$ is isomorphic to $\mathbb{Z}$ then the
$d$--fold cover $(B^4-\Delta)^d$ of $B^4-\Delta$ is a homology
circle.
\end{Lemma}

\begin{proof}
Let $(S^3-K)^d$ and $(B^4-\Delta)^d$ be the $d$--fold covers of
$(S^3-K)$ and $(B^4-\Delta)$ according to the following
homomorphisms $\varphi_1$, $\varphi_2$:
$$\bfig
  \barrsquare<800,400>[\pi_1(S^3-K)`H_1(S^3-K)`
    \pi_1(B^4-\Delta)`H_1(B^4-\Delta);
    \varphi_1`i_{*}`\cong`\varphi_2]
  \barrsquare(800,0)/->```->/<800,400>[\phantom{H_1(S^3-K)}`\mathbb{Z}/d`
    \phantom{H_1(B^4-\Delta)}`\mathbb{Z}/d;```]
  \morphism(0,400)|r|<0,-400>[\phantom{H_1(S^3-K)}`
  \phantom{H_1(B^4-\Delta)};\text{surj}]
  \efig$$
Since $K$ is a ribbon knot, $i_{*}\co \pi_1(S^3-K)\to\pi_1(B^4-\Delta)$
is surjective. It follows that the map $H_1((S^3-K)^d)\to
H_1((B^4-\Delta)^d)$ between the $d$--fold coverings is surjective
since $i_{*}(\ker\varphi_1)$ maps to the trivial element of
$\mathbb{Z}/d$ under $\varphi_2$.  Since $H_1((S^3-K)^d)$ is
isomorphic to $\mathbb{Z}$,  so is $H_1((B^4-\Delta)^d)$. To show
$H_{*}((B^4-\Delta)^d)=0$ for $*>1$, we consider the long exact
sequence of the pair $((B^4-\Delta)^d,\partial(B^4-\Delta)^d)$.
\begin{equation*}
\begin{split}
&H_4((B^4-\Delta)^d,\partial(B^4-
\Delta)^d)\stackrel{\partial_4}\rightarrow
H_3(\partial(B^4-\Delta)^d)\stackrel{i_3}\rightarrow
H_3(B^4-\Delta)^d \\
&\quad \stackrel{j_3}\rightarrow
H_3((B^4-\Delta)^d,\partial(B^4-\Delta)^d)
\stackrel{\partial_3}\rightarrow H_2(\partial(B^4-\Delta)^d)
\stackrel{i_2}\rightarrow H_2(B^4-\Delta)^d\rightarrow\cdots
\end{split}
\end{equation*}
Since $\partial_4$ is an isomorphism, $j_3$ is injective so that
$H_3((B^4-\Delta)^d)$ is isomorphic to $\im j_3=\ker \partial_3$.
Our claim is that $\partial_3\co H_3((B^4-\Delta)^d,\partial(B^4-\Delta)^d)
\rightarrow H_2(\partial(B^4-\Delta)^d)$ is
an isomorphism. Observe that
$\partial(B^4-\Delta)^d=(S^3-K)^d\cup\wwtilde \Delta\times\partial
D^2$ where $\wwtilde \Delta$ is the lifted disk of $\Delta$ in the
$d$--fold cover of $B^4$. By Poincar\'e Duality and the Universal
Coefficient Theorem,
$$H_3((B^4-\Delta)^d,\partial(B^4-\Delta)^d)\cong
H^1((B^4-\Delta)^d)\cong \Hom(H_1((B^4-\Delta)^d), \mathbb{Z})$$
and
\begin{align*}
H_2((S^3-K)^d\cup\wwtilde\Delta\times\partial D^2)&\cong
H^1((S^3-K)^d\cup\wwtilde\Delta\times\partial D^2)\\
&\cong  \Hom(H_1((S^3-K)^d\cup\wwtilde\Delta\times\partial
D^2),\mathbb{Z}).
\end{align*}
Since $H_1((B^4-\Delta)^d)$ and $H_1((S^3-K)^d)$ are isomorphic to
the group $\mathbb{Z}$ generated by the lifted meridian
$\wtilde\mu$ of $K$ in $S^3$,
$$H_3((B^4-\Delta)^d,\partial(B^4-\Delta)^d)\cong
H_2((S^3-K)^d\cup\wtilde\Delta\times\partial D^2)\cong\mathbb{Z}$$
and moreover the boundary map $\partial_3$ induced by the
restriction map from $(B^4-\Delta)^d$ to $(S^3-K)^d$. Hence
$\partial_3$ is an isomorphism and so this proves that
    $H_3((B^4-\Delta)^d)=0$ and also $H_4((B^4-\Delta)^d)=0$.

Considering that the Euler characteristic of $(B^4-\Delta)^d$ is
$\chi(B^4-\Delta)^d=d\cdot \chi(B^4-\Delta)$ and $H_*(S^3-K)\to
H_*(B^4-\Delta)$ is an isomorphism, we get $H_2((B^4-\Delta)^d)=0$.
\end{proof}

\begin{Rmk} We may look at \fullref{ex:4.6} to
see infinitely many knots whose $d$--fold covers satisfy the
condition in \fullref{lem:HC}.
\end{Rmk}

In the following Proposition, we will show that $W$ in
\fullref{def:cobo} is a homology cobordism. The condition
that $K$ is a ribbon knot
allows us to show that it is in fact a relative $h$--cobordism.

\begin{Prop} \label{prop:cobo}
If $K$ is a ribbon knot and the homology of $d$--fold cover
$(S^3-K)^d$ of $S^3-K$, $H_1((S^3-K)^d)\cong\mathbb{Z}$ with
$d\equiv\pm 1 \pmod{m}$ then there exists an $h$--cobordism $W$
between $M_0=X-\Sigma$ and $M_1=X-\Sigma_{K}(m)$ rel $\partial$.
\end{Prop}

\begin{proof}
Keeping the previous notation in mind, let's denote $W=X\times
I-(\Sigma\times I )_{A}(m)$, $M_0=X-\Sigma$ and $M_1=X-\Sigma_K(m)$.
To show that $W$ is $H_*$--cobordism rel $\partial$, we'll prove
$H_{*}(W,M_1)=H_{*}(W,M_0)=0$.

First, we need to describe $W$ and $M_1$ as follows; if we take a
neighborhood of the curve $\alpha$ in $\Sigma$ as $S^1\times B^3$
meeting $\Sigma$ on $S^1\times I$ then denoting the complement of
$S^1\times I$ in $\Sigma$ by $\Sigma_0$, we may write
\begin{equation}\label{cobo}
W=(X-S^1\times B^3-\Sigma_0)\times I\cup S^1\times_{\tau^m}
(B^3\times I-A_{+})
\end{equation}
and
\begin{equation}\label{cobom}
M_1=(X-S^1\times B^3-\Sigma_0)\cup S^1\times_{\tau^m_{K}}
(B^3-K_{+})
\end{equation}
Then considering the above description, the relative Mayer--Vietoris
sequence shows
$$H_{*}(W,M_1)\cong
 H_{*}(S^1\times_{\tau^m}(B^3\times I-A_{+}),
 S^1\times_{\tau_{K}^m}(B^3-K_{+})).$$
By the Alexander Duality, this relative homology group is same as
$$H_{*}(S^1\times_{\tau^m}(B^3\times
I,A_{+}),S^1\times_{\tau_{K}^m}(B^3,K_{+}))$$
which is trivial.
Similarly, we can show that $H_{*}(W,M_0)$ is trivial as well.

A similar argument shows that $H_{*}(V,\partial{M_0})$ is trivial
and hence we have shown that $W$ is a homology cobordism from $M_0$
to $M_1$ rel $\partial$. To assert that $W$ is a relative
$h$--cobordism, we need to show that $\pi_1(W)=\pi_1(X\times
I-\nu(\Sigma\times
I)_{A}(m))\cong\mathbb{Z}/d$.

For simplicity let us denote $U=X-S^1\times B^3-\Sigma_0$ and
$V=S^1\times_{\tau^m} (B^3-K_{+})$ in the decomposition
$(X-S^1\times B^3-\Sigma_0)\cup{S^1\times_{\tau^m} (B^3-K_{+})}$ of
$X-\Sigma_{K}(m)$. Then $U\cap{V}= S^1\times(\partial{B^3}-\{
\mbox{two points} \}) $. Denoting $V'= S^1\times_{\tau^m}(B^3\times
I- A_{+})$, we also rewrite
$$W=(X-S^1\times B^3-\Sigma_0)\times I \cup
S^1\times_{\tau^m}(B^3\times I-A_{+})=U\times I \cup V'.$$
Then the intersection $U\times I \cap V'$ is $S^1\times(\partial{B^3}-\{$
two points $\})\times I=(U\cap V)\times I $.

Applying Van Kampen's theorem for these decompositions of $M_1$ and
$W$, we have the two commutative diagrams:
$$\bfig
\barrsquare<1200,500>[\pi_1 (U\cap V)`\pi_1(U)`
  \pi_1(S^1\times_{\tau^m} (B^3-K_{+}))`\pi_1(X-\Sigma_{K}(m));
  \varphi_1`\varphi_2`\psi_2`\psi_2]
  \efig$$
and
$$\bfig
\barrsquare<1500,500>[\pi_1 ((U\cap V)\times I)`\pi_1(U \times I)`
  \pi_1(S^1\times_{\tau^m} (B^3\times I -A_+))`
  \pi_1(X\times I - (\Sigma\times I)_A(m));
  \varphi_1'`\varphi_2'`\psi_2'`\psi_2']
  \efig$$
Let
\begin{align*}
i_{1}&\co \pi_1 (U\cap V)\rightarrow
  \pi_1 ((U\cap V)\times I)\\
i_{2}&\co \pi_1 (U)\rightarrow
  \pi_1 (U\times I)\\
i_3&\co \pi_1(S^1\times_{\tau^m} (B^3-K_{+}))\rightarrow
\pi_1(S^1\times_{\tau^m} (B^3\times I-A_{+}))
\end{align*}
be the maps induced
by inclusions. Then clearly $i_{1}$ and $i_{2}$ are isomorphisms. To
show that $i_3$ is surjective, let's consider the fundamental group
of mapping cylinders $S^1\times_{\tau^m}(B^3-K_{+})$ and
$S^1\times_{\tau^m} (B^3\times I-A_{+})$. Then representing the
element $[S^1]$ in the fundamental group as $\alpha'$, we present
\begin{equation*}
\begin{split}
\pi_1(S^1\times_{\tau^m}&(B^3-K_{+}))=\\
&\langle\pi_1(B^3-K_+), \alpha' ~|~
\alpha'^{-1}\beta\alpha' =\tau_{K*}^{m}(\beta)\text{ for all } \beta \in
\pi_1(B^3-K_+)\rangle.
\end{split}
\end{equation*}
and
\begin{equation*}
\begin{split}
\pi_1(S^1 & \times_{\tau^m}(B^3\times I-A_{+}))=\\
 & \langle\pi_1(B^3\times
I-A_{+}), \alpha' ~|~ \alpha'^{-1}\beta'\alpha'
=\tau_{*}^{m}(\beta') \text{ for all }
\beta' \in \pi_1(B^3\times I-A_{+})\rangle.
\end{split}
\end{equation*}
Since $K$ is a ribbon knot,
$\pi_1(S^3-K)\longrightarrow\pi_1(S^3\times I-A)$ is surjective. So
is $i_3$. Then by chasing the diagram, we have a surjective map
$$\pi_1(X-\Sigma_{K}(m)) \rightarrow \pi_1(X\times I-(\Sigma\times
I)_{A}(m)).$$
By \fullref{prop:FG},
$\pi_1(X-\Sigma_{K}(m))=\mathbb{Z}/d$. Since $W$ is an
$H_*$--cobordism by the above argument, $H_1(X\times I-(\Sigma\times
I)_{A}(m))=\mathbb{Z}/d$ so that $\pi_1(X\times I-(\Sigma\times
I)_{A}(m))=\mathbb{Z}/d$.

Now let us  prove that the inclusion $i\co M_1\longrightarrow W$ is a
homotopy equivalence. The above work shows that the induced map
$i_{*}\co \pi_1{M_1}\longrightarrow \pi_1{W}\cong\mathbb{Z}/d$ is an
isomorphism. So, the $d$--fold covers ${W}^d$ and ${M_1}^d$ of $W$ and
$M_1$ become universal covers and so we denote $\wwtilde W={W}^d$,
$\wwtilde M_1={M_1}^d$. Then we claim that the inclusion $\wwtilde
M_1\to \wwtilde W$ induces an isomorphism in homology.
Considering the decompositions of $W$ and $M_1$ in \eqref{cobo} and
\eqref{cobom}, we can express their $d$--fold covers as the $d$--fold
covers of subcomponents associated to their inclusion maps to
$H_1(W)\cong\mathbb{Z}/d$:
\begin{equation*}
\begin{split}
\wwtilde W&=(X\times I-(\Sigma\times I)_{A}(m))^d\\
&=((X-S^1\times B^3-\Sigma_0)\times I)^d\cup
(S^1\times_{\tau^{m}}(B^3\times I-A_{+}))^d
\end{split}
\end{equation*}
and
$$\wwtilde M_1=(X-\Sigma_{K}(m)))^d= (X-S^1\times B^3-\Sigma_0)^d \cup
(S^1\times_{\tau_{K}^{m}}(B^3-K_{+}))^d.$$
In the inclusion-induced map $j\co H_1(S^1\times_{\tau^m}(B^3\times I-
A_{+}))\longrightarrow H_1(W)\cong\mathbb{Z}/d$, from our choice of
the curve $\alpha$ in $\Sigma$ mentioned in the beginning of
\fullref{sec:3}, we can easily check that in the Mayer--Vietoris sequence, the
homology element $[S^1\times pt\times 0]$ with $pt\in(\partial
B^3-\text{two points})$ maps under $j$ to a trivial element in $H_1(W)$.
The Mayer--Vietoris sequence for the decomposition of $W$ follows.
\begin{equation*}
\begin{split}
\cdots &\longrightarrow H_1(S^1\times(\partial B^3-\{\text{two
points}\})\times I)
\stackrel{\varphi}\longrightarrow \\
&\stackrel{\varphi}\longrightarrow H_1((X-S^1\times
B^3-\Sigma_0)\times I)\oplus H_1(S^1\times_{\tau^m}(B^3\times I-
A_{+}))
\stackrel{\psi}\longrightarrow \\
& \stackrel{\psi}\longrightarrow H_1(W)\longrightarrow 0 .
\end{split}
\end{equation*}
First we note the image of a generator $[S^1\times pt\times
0]\in H_1(S^1\times(\partial B^3-\{\text{two points}\})\times I)$ under
$\varphi$ is  $(0,[S^1\times pt\times 0]) \in H_1((X-S^1\times
B^3-\Sigma_0)\times I)\oplus H_1(S^1\times_{\tau^m}(B^3\times I-
A_{+})) $ since the pushed-off curve of $\alpha$ along a
trivialization is zero in $H_1((X-S^1\times B^3-\Sigma_0)\cong
H_1(X-\Sigma)$ by adjusting the framing of the curve $\alpha$.

So, since $(0,[S^1\times pt\times 0])$ is in the kernel of $\psi$,
$[S^1\times pt\times 0]$ maps to the trivial element in $H_1(W)$.
Then we know the $d$--fold cover of $S^1\times_{\tau^{m}}(B^3\times
I-A_{+})$ has the form
$$S^1\times_{\wtilde{\tau}^{m}}(B^3\times I-A_{+})^d$$
for a proper lifted map $\wtilde{\tau}^{m}$ of $\tau^{m}$ and by the
same reason, the $d$--fold cover of
$$S^1\times_{\tau_{K}^{m}}(B^3-K_{+})$$
is also of the form
$$S^1\times_{\wtilde{\tau}_{K}^{m}}(B^3-K_{+})^d$$
for some lift $\wtilde{\tau}_{K}^{m}$ of $\tau_{K}^{m}$.

Then we have a simple form of the relative homology of the pair
$(W,M_1)$,
$$H_{*}(\wwtilde{W},\wwtilde{M}_1)\cong
H_{*}(S^1\times_{\wwtilde\tau^{m}}(B^3\times I-A_{+})^d,
S^1\times_{\wwtilde\tau_{K}^{m}}(B^3-K_{+})^d).$$
Since $K$ is a ribbon knot and $H_1((S^3-K)^d)\cong\mathbb{Z}$, it follows
by \fullref{lem:HC} that
$H_*\bigl((B^3{\times}I{-}A_+)^d, (B^3{-}K_+)^d\bigr)=0$.
So,
the homology $H_{*}(\wwtilde{W},\wwtilde{M}_1)$ is trivial. By the
Whitehead theorem, we get
$\pi_{n}\wwtilde{M}_1\cong \pi_{n}\wwtilde{W}$
for $n>1$. Since $\pi_{n}\wwtilde{M}_1\cong\pi_{n}{M_1}$ and
$\pi_{n}\wwtilde{W}\cong\pi_{n}{W}$, it follows that
$i_{*}\co \pi_{n}M_1\to \pi_{n}W$ is an
isomorphism. Therefore, again by Whitehead's theorem,
$i\co M_1\longrightarrow W$ is a homotopy equivalence.
\end{proof}

Now we need to recall the definition of torsion, as given in \cite{M}
or \cite{T} to show the Whitehead torsion of the pair $(W,M_0)$
constructed above is zero.

Let $\Lambda$ be an associative ring with unit such that for any
$r\neq s\in \mathbb{N}$, $\Lambda^{r}$ and $\Lambda^{s}$ are not
isomorphic as $\Lambda$--modules. Consider an acyclic chain complex
$C$ of length $m$ over $\Lambda$ whose chain groups are finite free
$\Lambda$--modules with a preferred basis $c_i$ for each chain
complex $C_i$. Then the torsion of the chain complex $C$ --- written
$\tau(C)$ --- is defined as follows.

Let $GL(\Lambda)=\bigcup_{n\geq 0} GL(n,\Lambda)$ be the infinite
general group. The torsion $\tau(C)$ will be an element of the
abelianization of $GL(\Lambda)$, denoted by $K_1(\Lambda)$.
  Pick ordered bases $b_i$
of $B_i=\Im\partial_{i}$ and combine them to bases $b_ib_{i-1}$ of
$C_i$. For the distinguished basis $c_i$ of $C_i$, let
$(b_ib_{i-1}/c_i)$ is the transition matrix over $\Lambda$.
Denoting the corresponding element of $K_1(\Lambda)$ by
$[b_ib_{i-1}/c_i]$, define the torsion
$$\tau(C)=\prod_{i=0}^{m}[b_ib_{i-1}/c_i]^{(-1)^{i+1}}\in
K_1(\Lambda).$$
In particular, if $(K,L)$ is a pair of finite, connected CW
complexes such that $L$ is a deformation retract of $K$ then
consider the universal covering complexes $\wwtilde K\supset\wtilde
L$ of $K$ and $L$. Let's denote $\pi$ by the fundamental group of
$K$. Then we obtain an acyclic free chain
$\mathbb{Z}[\pi]$--complex $C(\wwtilde K,\wtilde L)$. So we have a
well defined torsion $\tau=\tau(K,L)$ in the Whitehead group
$\Wh(\pi)=K_1(\mathbb{Z}[\pi])/\pm\pi$, the so-called
`Whitehead torsion'.

The $h$--cobordism we have constructed is built out of several
pieces, and so our strategy is to compute the Whitehead torsion in
terms of those pieces. The pieces may not be $h$--cobordisms, so they
don't have a well-defined Whitehead torsion. However, they do have a
more general kind of torsion, the Reidemeister--Franz torsion, which
we briefly outline. It will turn out that the Reidemeister--Franz
torsion of the pieces determines the Whitehead torsion of the
$h$--cobordism. Moreover, the Reidemeister--Franz torsion satisfies
gluing laws which will be able us to compute its value in terms of
the pieces.

 The `Reidemeister--Franz torsion' is defined as follows. Consider
the pair $(K,L)$ of finite, connected CW-complexes but not requiring
that $L$ is a deformation retract of $K$. Then keeping the notation
above,  the cellular chain group $C_i(\wwtilde K,\wtilde L)$ is a free
$\mathbb{Z}[\pi]$--module as before. Let $\Lambda$ be an associative
ring with unit with the above property. Given a ring homomorphism
$\varphi\co \mathbb{Z}[\pi]\longrightarrow\Lambda$, consider a free
chain complex
$$C^{\varphi}(K,L)=\Lambda\otimes_{\varphi}C(\wwtilde K,\wtilde L).$$
If $C^{\varphi}$ is acyclic, the torsion corresponding the chain
complex $C^{\varphi}$ is well defined. We will denote
$\tau^{\varphi}(K,L)\in K_1(\Lambda)/\pm\varphi(\pi)$. If $\Lambda$
is a field then $K_1(\Lambda)=\Lambda^*$ so that
$\tau^{\varphi}(K,L)\in \Lambda^*/\pm\varphi(\pi)$.

If the original complex $C$ is acyclic then the new complex
$C^\varphi$ is also acyclic and so when the Whitehead torsion of
$(K,L)$ is defined, the Reidemeister torsion of $(K,L)$ is also
defined associated to the identity homomorphism
$\id\co \mathbb{Z}[\pi]\longrightarrow \mathbb{Z}[\pi]$. However the
relation
$$\tau^{\varphi}(K,L)=\varphi_*\tau(K,L)$$
shows that if the Reidemeister torsion associated to the identity is
trivial then the Whitehead torsion is zero. We also need to know
some formulas to compute torsion. Suppose $K=K_1\cup K_2$,
$K_0=K_1\cap K_2$, $L=L_1\cup L_2$, $L_0=L_1\cap L_2$ and that
$i\co L\longrightarrow K$ is the inclusion which is restricted to
homotopy equivalences $i_{\alpha}\co L_{\alpha}\longrightarrow
K_{\alpha}$  (for $\alpha=0,1,2$). Then $i$ is a homotopy equivalence and
we have a formula called the `sum theorem' in Whitehead torsion (see
\cite{T})
$$\tau(K,L)=i_{1*}\tau(K_1,L_1)+i_{2*}\tau(K_2,L_2)-i_{0*}\tau(K_0,L_0).$$
Using the multiplicativity of the torsion and the Mayer--Vietoris
sequence we obtain a similar one called the `gluing formula' in the
Reidemeister torsion (see \cite{T}).

Given subcomplexes $X_1$ and $X_2$ of $X$ such that $X=X_1\cup X_2$ and
$X_1\cap X_2=Y$, let $\varphi\co \mathbb{Z}[H_1(X)]\longrightarrow
\Lambda$ be a ring morphism where $\Lambda$ is a ring as above. Let
$i\co \mathbb{Z}[H_1(Y)]\longrightarrow \mathbb{Z}[H_1(X)]$ and
$i_{\alpha}\co \mathbb{Z}[H_1(X_{\alpha})]\longrightarrow
\mathbb{Z}[H_1(X)]$  (for $\alpha=1,2$) denote the inclusion-induced
morphisms. If $\tau^{\varphi\circ i}(Y)\neq 0$ then we have the
gluing formula
$$\tau^{\varphi}(X)\cdot\tau^{\varphi\circ i}(Y)= \tau^{\varphi\circ
i_1}(X_1)\cdot \tau^{\varphi\circ i_2}(X_2).$$
Now considering our situation, we have shown that $W$ is a relative
$h$--cobordism from $M_0$ to $M_1$ with $\pi_1(W)\cong \mathbb{Z}/d$
and so the Whitehead torsion $\tau(W,M_0)\in Wh(\mathbb{Z}/d)$ is
defined. Recall that the decomposition of the pair
$$(W,M_0)=(X\times I-(\Sigma\times I )_{A}(m),X-\Sigma)$$
in \eqref{cobo} and \eqref{cobom} is
$$((X{-}S^1\times
B^3{-}\Sigma_0)\times I\cup S^1\times_{\tau^m}(B^3\times I{-} A_{+}),
X{-}S^1\times B^3{-}\Sigma_0\cup S^1\times (B^3{-}I)).$$
If we rewrite this as
$$((X{-}S^1{\times} B^3{-}\Sigma_0){\times} I,
  X{-}S^1{\times} B^3{-}\Sigma_0) \cup
 (S^1{\times_{\tau^m}}(B^3{\times} I{-} A_{+}), S^1{\times} (B^3{-}I)),$$
then we can observe that the Whitehead torsion of the first
component pair
$$((X-S^1\times B^3-\Sigma_0)\times I, X-S^1\times B^3-\Sigma_0)$$
is zero and so we would like to
attempt to use the sum theorem for this decomposition. But in the
second pair, $S^1\times_{\tau^m}(B^3\times I- A_{+})$ is just a
homology cobordism which means $S^1\times (B^3-I)$ may not be a
deformation retract of $S^1\times_{\tau^m}(B^3\times I- A_{+})$.
Then the Whitehead torsion $\tau(S^1\times_{\tau^m}(B^3\times I-
A_{+}), S^1\times (B^3-I))$ is not defined and thus we can not apply
the sum theorem in order to show the Whitehead torsion
$\tau(W,M_0)=0$. But we will show later that
$\tau(S^1\times_{\tau^m}(B^3\times I- A_{+}), S^1\times (B^3-I))$
is well defined under an additional assumption to make the complex of
the $d$--fold cover of the pair, $C((S^1\times_{\tau^m}(B^3\times I-
A_{+}))^d, (S^1\times (B^3-I))^d)$, acyclic with
$\mathbb{Z}[\mathbb{Z}/d]$ coefficient. So instead of computing the
Whitehead torsion, we will show that the Reidemeister torsion
$\tau^{\id}(W,M_0)$, denoted simply by $\tau(W,M_0)$, according to
  the coefficient
$\mathbb{Z}[\mathbb{Z}/d]$ is trivial. Applying the gluing formula
to the above decomposition instead of the sum theorem, we can obtain
a simpler method to compute the Reidemeister torsion for the pair
$(W,M_0)$.

Now we first need to consider the torsion of certain fibration
over a
circle with a homologically trivial fiber.

A relative fiber bundle
$$(F,F_0)\hookrightarrow (X,Y)\stackrel{\pi}\longrightarrow S^1$$
means that $F\hookrightarrow X\stackrel{\pi}\longrightarrow S^1$ is
a fiber bundle with a trivialization $\{\varphi_{\alpha},
U_{\alpha}\}$ satisfying that for an open cover $U_\alpha\subset
S^1$,  $(\pi^{-1}(U), Y\cap \pi^{-1}(U))\cong U\times (F,F_0)$ and
the diagram
$$\bfig
  \Vtriangle<500,400>[(\pi^{-1}(U), Y\cap \pi^{-1}(U))`
    U{\times}(F,F_0)`U;
    \varphi_\alpha``]
  \efig$$
commutes. We will now prove the following result.

\begin{Prop} \label{prop:torsion}
Let $(F,F_0)\hookrightarrow (X,Y)\longrightarrow S^1$ be a smooth,
relative fiber bundle over $S^1$ such that the fiber pair $(F,F_0)$
is homologically trivial. Suppose that $G$ is a group and $\rho \co 
H_1(X)\longrightarrow G$ is a group homomorphism such that the image
under $\rho$ of the homology class $[S^1]\in H_1(X)$ of the base
space in the fibration has finite order in $G$. Let $(\wwtilde
F,\wwtilde F_0)$ be the cover of $(F,F_0)$ associated to the
homomorphism
$$H_1(F)\hookrightarrow
H_1(X)\stackrel{\rho}\longrightarrow G$$ and denote again by $\rho$
the induced map $\mathbb{Z}[H_1(X)]\longrightarrow \mathbb{Z}[G]$.
If the cover $(\wwtilde F,\wwtilde F_0)$ is homologically trivial, that
is $H_*(F,F_0; \mathbb{Z}[G])=0$ then the torsion
$\tau^{\rho}(X,Y)\in K_1(\mathbb{Z}[G])/\pm G$ is trivial.
\end{Prop}

\begin{proof}
We may assume that $X$ is a mapping torus $X=S^1\times_{\varphi}F$
with the monodromy map $\varphi$ of the fibration. Let $(\wwtilde
X,\wwtilde Y)$ be the cover of $(X,Y)$ associated to $\rho$. Then
$\wwtilde X$ is also a mapping torus since the homology image $\rho
([S^1])$ is of finite order in $G$. So, let us say $\wwtilde
X=S^1\times_{\wwtilde \varphi}\wwtilde F$ where $\wwtilde F$ is the cover
associated to $H_1(F)\hookrightarrow
H_1(X)\stackrel{\rho}\longrightarrow  G$ and $\wtilde\varphi$ is a
lift of $\varphi$ in $\wwtilde X$. Similarly, we also say $\wwtilde Y=
S^1\times_{\wtilde \varphi}\wwtilde F_{0}$. Considering the Wang exact
homology sequence and the  Five Lemma we have
$$\begin{array}{ccccccccc}
H_*(\wwtilde F_0) \! &
\stackrel{\wtilde\varphi_{*}^{-1}}\longrightarrow & \! H_*(\wwtilde
F_0) \! & \longrightarrow & \! H_*(S^1\times_{\wtilde
\varphi}\wwtilde F_0) \! & \longrightarrow & \! H_{*-1}(\wwtilde F_0)
\! &\stackrel{\wtilde\varphi_{*}^{-1}}\longrightarrow
& \! H_{*-1}(\wwtilde F_0)\vspace{5 pt} \\
\cong\!\big\downarrow & & \cong\!\big\downarrow & & \big\downarrow
& & \cong\!\big\downarrow &
& \cong\!\big\downarrow \\
H_*(\wwtilde F) \! &
\stackrel{\wtilde\varphi_{*}^{-1}}\longrightarrow & \! H_*(\wwtilde
F) \! & \longrightarrow & \! H_*(S^1\times_{\wtilde \varphi}\wwtilde
F) \! & \longrightarrow & \! H_{*-1}(\wwtilde F) \!
&\stackrel{\wtilde\varphi_{*}^{-1}}\longrightarrow
& \! H_{*-1}(\wwtilde F)\\
\end{array}$$
and we get an acyclic complex  $C_*(S^1\times_{\wtilde \varphi}\wwtilde
F,S^1\times_{\wtilde \varphi}\wwtilde F_0)$ since $H_{*}(\wwtilde
F_{0})\longrightarrow H_{*}(\wwtilde F)$ is an isomorphism. Thus, the
associated torsion $\tau^{\rho}(X,Y)$ is defined.

Now we consider the Mayer--Vietoris sequence for $(\wwtilde
X,\wwtilde Y)=(S^1\times_{\wtilde\varphi}\wwtilde F,
S^1\times_{\wtilde\varphi}\wwtilde F_{0})$. Let us consider closed
manifold pairs $(X_1, Y_1)=([0,\frac{1}{2}]\times \wwtilde F,
[0,\frac{1}{2}]\times \wwtilde F_{0})$ and $(X_2, Y_2)=([\frac{1}{2},
1]\times \wwtilde F, [\frac{1}{2}, 1]\times \wwtilde F_{0})$. Define a
map $f$ of a subspace $A:=\{0\}\times\wwtilde
F\cup\{\frac{1}{2}\}\times\wwtilde F$ of $X_1$ into $X_2$ by
$$f|_{\{0\}\times\wwtilde F}=\wtilde\varphi\times\{1\},
f|_{\{1/2\}\times\wwtilde F}=1_{\{1/2\}\times\wwtilde F}.$$
Then letting  $B:=\{0\}\times\wwtilde
F_{0}\cup\{\frac{1}{2}\}\times\wwtilde F_{0}\subset A$, we can
consider $(S^1\times_{\wtilde \varphi}\wwtilde F, S^1\times_{\wtilde
\varphi}\wwtilde F_{0})$ as the adjunction space $(X_1\cup_{f}X_2,
Y_1\cup_{f}Y_2)$ of the system $(X_1, Y_1)\supset
(A,B)\stackrel{f}\longrightarrow (X_2, Y_2)$. There is a short
exact sequence
\begin{multline*}
0  \longrightarrow C_*(X_1\cap X_2,Y_1\cap Y_2)\longrightarrow
C_*(X_1, Y_1)\oplus C_*(X_2, Y_2) \\
 \longrightarrow C_*(X_1\cup_{f}X_2,
Y_1\cup_{f}Y_2)\longrightarrow 0.
\end{multline*}
If we rewrite this then we have
\begin{multline*}
0\longrightarrow C_*(\wwtilde F, \wwtilde F_{0})\oplus C_*(\wwtilde
F, \wwtilde F_{0})) \\
\longrightarrow C_*([0,1/2]\times \wwtilde F, [0,1/2]\times \wwtilde
F_0)\oplus C_*([1/2,1]\times \wwtilde F, [1/2,1]\times
\wwtilde F_0)\\
\longrightarrow C_*(S^1\times_{\wtilde \varphi}\wwtilde F,
S^1\times_{\wtilde \varphi}\wwtilde F_0)\longrightarrow 0.
\end{multline*}
If $(\wwtilde F,\wwtilde F_0)$ is homologically
trivial, it follows that if $j\co \mathbb{Z}[H_1(F)]\longrightarrow
\mathbb{Z}[H_1(X)]$ denotes the morphism induced by inclusion then
the torsion $\tau^{\rho\circ j}(F, F_0)$ is defined. From the above
short exact sequence and the multiplicativity of the torsion we
deduce that
$$\tau^{\rho\circ j}(F, F_0)\cdot \tau^{\rho\circ j}(F,F_0)=
(\tau^{\rho\circ j}(F, F_0)\cdot \tau^{\rho\circ
j}(F,F_0))\cdot\tau^{\rho}
(S^1\times_{\varphi}F,S^1\times_{\varphi}F_0).$$
This implies that $\tau^{\rho}
(S^1\times_{\varphi}F,S^1\times_{\varphi}F_0)=\tau^{\rho}(X,Y)\in
K_1(\mathbb{Z}[G])/\pm G$ is trivial.
\end{proof}

Using the proposition above, we get topological equivalence
classes of $(X,\Sigma_K(m))$ under the following condition.

\begin{Thm} \label{thm:last}
If $K$ is a ribbon knot and the homology of $d$--fold cover
$(S^3-K)^d$ of $S^3-K$, $H_1((S^3-K)^d)\cong\mathbb{Z}$ with
$d\equiv\pm 1 \pmod{m}$ then $(X,\Sigma)$ is pairwise homeomorphic
to $(X,\Sigma_{K}(m))$.
\end{Thm}

\begin{proof}
Under these assumptions, we have a relative $h$--cobordism $W$ from
$M_0=X-\Sigma$ to $M_1=X-\Sigma_{K}(m)$ by
\fullref{prop:cobo}.  As we discussed before, in order to
show the Whitehead torsion $\tau(W,M_0)=0\in  Wh(\mathbb{Z}/d)$, it
is sufficient to show that the Reidemeister torsion $\tau(W,M_0)\in
Wh(\mathbb{Z}/d)$ associated to the identity map $\id\co 
\mathbb{Z}[\mathbb{Z}/d]\longrightarrow
\mathbb{Z}[\mathbb{Z}/d]$ is trivial.

Consider the decomposition of the pair $(W,M_0)$,
$$((X-S^1\times B^3-\Sigma_0)\times I,
  X-S^1\times B^3-\Sigma_0)\cup
(S^1\times_{\tau^m}(B^3\times I- A_{+}), S^1\times (B^3-I)).$$
To apply the gluing formula of the Reidemeister torsion for this
decomposition, we need to check the torsion of each component is
defined.

First, the torsion $\tau((X-S^1\times B^3-\Sigma_0)\times I,
  X-S^1\times B^3-\Sigma_0)$ is clearly defined and trivial.
To check the torsion of the second component, we will show the
relative chain complex $C((S^1\times_{\tau^m}(B^3\times I-
A_{+}))^d, (S^1\times
(B^3-I))^d)$ of $d$--fold covers is acyclic.

The same argument in the proof of \fullref{prop:cobo} shows
that the $d$--fold cover
$$(S^1\times_{\tau^m}(B^3\times I- A_{+}))^d$$
associated to the inclusion-induced map
$$j\co H_1(S^1\times_{\tau^m}(B^3\times I- A_{+}))\longrightarrow
H_1(W)\cong\mathbb{Z}/d$$
is a mapping torus with the $d$--fold cover of
$B^3\times I- A_{+}$ that is $S^1\times_{\wtilde\tau^m}(B^3\times I-
A_{+})^d$. Similarly, the $d$--fold cover
$(S^1\times_{\tau^m}(B^3-I))^d$ is
$S^1\times_{\wtilde\tau^m}(B^3-I)^d$.

Observing the proof of \fullref{lem:HC}, we have an isomorphism
$H_*((S^3-K)^d)\longrightarrow H_*((B^4-\Delta)^d)$ when $K$ is a
ribbon knot and $H_1((S^3-K)^d)\cong\mathbb{Z}$. In other words,
$H_*((B^4-\Delta)^d, (S^3-K)^d)=0$. Excision argument shows that
this is isomorphic to
$$H_*((B^3\times I-A_+)^d, (B^3-K_+)^d)=0\cong H_*((B^3\times I-A_+)^d,
(B^3-I)^d).$$
This gives that
\begin{equation*}
\begin{split}
H_*((&S^1\times_{\tau^m} (B^3\times
I- A_{+}))^d,(S^1\times_{\tau^m}(B^3-I))^d) \\
&=\ H_*(S^1\times_{\wtilde\tau^m}(B^3\times
I-A_+)^d,S^1\times_{\wtilde\tau^m} (B^3-I)^d)=0.
\end{split}
\end{equation*}
Then the torsion $\tau^j(S^1\times_{\tau^m}(B^3\times I-
A_{+}),S^1\times_{\tau^m}(B^3-I))$ associated to the induced ring
homomorphism $j\co \mathbb{Z}[H_1(S^1\times_{\tau^m}(B^3\times I-
A_{+}))]\longrightarrow
\mathbb{Z}[H_1(W)]\cong\mathbb{Z}[\mathbb{Z}/d]$
is defined.

Now applying the gluing formula of the Reidemeister torsion for the
decomposition of $(W,M_0)$, we have
\begin{equation*}
\begin{split}
\tau(W&,M_0) \cdot\tau(\partial(X-S^1\times B^3-\Sigma_0)\times I,
  \partial(X-S^1\times B^3-\Sigma_0))\\
&=\ \tau((X-S^1\times B^3-\Sigma_0)\times I,
X-S^1\times B^3-\Sigma_0)\ \cdot \\
&\qquad\qquad\qquad \tau(S^1\times_{\tau^m}(B^3\times I- A_{+}),
S^1\times (B^3-I)).
\end{split}
\end{equation*}
Hence,
$$\tau(W,M_0)=\tau(S^1\times_{\tau^m}(B^3\times I- A_{+}), S^1\times
(B^3-I)).$$
To compute $\tau(S^1\times_{\tau^m}(B^3\times I- A_{+}), S^1\times
(B^3-I))$, we note that
$$(B^3\times I- A_{+}, B^3-I)\hookrightarrow
(S^1\times_{\tau^m}(B^3\times I- A_{+}), S^1\times
(B^3-I))\longrightarrow S^1$$
is a smooth fiber bundle over $S^1$ with the fiber $(B^3\times I-
A_{+}, B^3-I)$. Clearly the fiber $(B^3\times I- A_{+}, B^3-I)$ is
homologically trivial and by the above argument, the $d$--fold cover
$((B^3\times I- A_{+})^d, (B^3-I)^d)$ associated to $j$ is also
homologically trivial.  Thus, by \fullref{prop:torsion} the
torsion $\tau(S^1\times_{\tau^m}(B^3\times I- A_{+}), S^1\times
(B^3-I))$ is trivial and thus the Whitehead torsion $\tau(W,M_0)=0$.
Then by Freedman's work \cite{Freedman}, the $h$--cobordism $W$ is
topologically trivial and so  the complements $X-\Sigma$ and
$X-\Sigma_K(m)$ are homeomorphic. The homeomorphism
$\partial\nu(\Sigma)\longrightarrow
\partial\nu(\Sigma_K(m))$ extends to a homeomorphism
$(X,\Sigma)\longrightarrow(X,\Sigma_K(m))$.
\end{proof}

\begin{Ex}
\label{ex:4.6}
Let's consider examples $(X, \Sigma_K(m))$ which are smoothly knotted
but top\-o\-log\-i\-cal\-ly standard.
Let $J$ be a torus knot $T_{p,q}$ in $S^3$ such that $p$ and $q$ are
coprime positive integers. Then we have a ribbon knot $K=J\# {-J}$
with its Alexander polynomial $\Delta_{K}(t)=(\Delta_{J}(t))^2$
where
$$\Delta_{J}(t)=\frac{(1-t)(1-t^{pq})}{(1-t^p)(1-t^q)}.$$
Note that the $d$--fold cover of $S^3$ branched along the torus knot
$J=T_{p,q}$ is the Brieskorn manifold $\Sigma(p,q,d)$, and that this
manifold is a homology sphere if p,q and  d are pairwise relatively
prime. Since $(S^3,K)^d$ is $\Sigma(p,q,d)\# \Sigma(p,q,d)$,
$(S^3,K)^d$ is an integral homology 3--sphere. We might obtain a
direct proof for this by computing the order of $H_1((S^3,K)^d)$ of
$d$--fold cover $(S^3,K)^d$ of $S^3$ branched over $K$. In fact, Fox
\cite{Fox} proved that
$$|H_1((S^3,K)^d)|=\prod_{i=0}^{d-1}\Delta_{K}(\zeta^i)$$
where $\zeta$ is a primitive $d$th root of unity. And it's easy to
show that
$$\prod_{i=0}^{d-1}\Delta_{K}(\zeta^i)=1.$$
So, we obtain a ribbon knot $K$ with $\Delta_K(t)\neq 1$ and the
$d$--fold branch cover $(S^3,K)^d$ is a homology 3--sphere when
$(p,d)=1$ and $(q,d)=1$. Then by \fullref{thm:main} and
\fullref{thm:last}, we have infinitely many pairs
$(X,\Sigma_K(m))$ which are smoothly knotted but not topologically.
\end{Ex}

\subsection*{Acknowledgements}

I thank the referee for pointing out a gap in the proof of one of
main theorems and other helpful comments. I would also like to thank
Fintushel and Stern for their note and I would like to express my
sincere gratitude to my advisor Daniel Ruberman for his tremendous
help and support.

\bibliographystyle{gtart}
\bibliography{link}

\end{document}